\documentclass[11pt]{article}
\usepackage{amsfonts}
\usepackage{amssymb}
\usepackage{amsmath}
\usepackage{amsxtra}
\usepackage{graphicx}
\usepackage{psfrag} 
\usepackage{bm}
\usepackage{mathrsfs}
\usepackage{stmaryrd}
\usepackage{subfig} 
\usepackage{tikz} 
\usepackage{empheq}
\usepackage[utf8]{inputenc}
\usepackage[english]{babel}
\usepackage[square,numbers]{natbib} 
\usepackage{pgfplots}   
\usepackage{pgfplotstable}
\pgfplotsset{compat=1.3}
\pgfplotsset{every axis legend/.append style={
    at={(1.05,1)},
    anchor=north west,font=\small}}
\usepackage{filecontents}
\usepackage{float} 
\usepgfplotslibrary{units}

\usepackage{caption} 
\usepackage{hyperref}  
\usepackage{hypcap} 
\usepackage{etoolbox}

\usepackage{xcolor}
  
\makeatletter
\patchcmd{\maketitle}{\@fnsymbol}{\@alph}{}{}  
\makeatother

\title{Modeling of dendritic solidification and numerical analysis of the phase-field approach to model complex morphologies in alloys}
\author{Kunal Bhagat \thanks{Department of Mechanical Engineering, University of Wisconsin-Madison} \and Shiva Rudraraju \thanks{Department of Mechanical Engineering, University of Wisconsin-Madison; corresponding author} \\[3ex]}
\date{Preprint}
\begin{document}
\maketitle

\section*{Abstract} 
Dendrites are one of the most widely observed patterns in nature and occur across a wide spectrum of physical phenomena. In solidification and growth patterns in metals and crystals, the multi-level branching structures of dendrites pose a modeling challenge, and a full resolution of these structures is computationally demanding. In the literature, theoretical models of dendritic formation and evolution, essentially as extensions of the classical moving boundary Stefan problem exist. Much of this understanding is from the analysis of dendrites occurring during the solidification of metallic alloys. Motivated by the problem of modeling microstructure evolution from liquid melts of pure metals and alloys during MAM, we developed a comprehensive numerical framework for modeling a large variety of dendritic structures that are relevant to metal solidification. In this work, we present a numerical framework encompassing the modeling of Stefan problem formulations relevant to dendritic evolution using a phase-field approach and a finite element method implementation. Using this framework, we model numerous complex dendritic morphologies that are physically relevant to the solidification of pure melts and binary alloys. The distinguishing aspects of this work are - a unified treatment of both pure metals and alloys; novel numerical error estimates of dendritic tip velocity; and the convergence of error for the primal fields of temperature and the order parameter with respect to numerical discretization. To the best of our knowledge, this is a first-of-its-kind study of numerical convergence of the phase-field equations of dendritic growth in a finite element method setting. Further, we modeled various types of physically relevant dendritic solidification patterns in 2D and 3D computational domains.

\section{Introduction}
Dendrites are tree-like patterns with complex multi-level branch structures that are observed across a wide spectrum of physical phenomena - from snow flakes to river basins; from bacterial colonies to lungs and vascular systems; and ubiquitously in solidification and growth patterns in metals and crystals. In the context of solidification problems involving a pure metal or metallic alloys, dendritic ``trees'' with primary and secondary branches grow from a nucleation point or a field perturbation. These solidification dendrites often initiate at random nucleation points or at domain boundaries at the initial stage, followed by anisotropic interface growth during the intermediate stage, eventually leading to grain formation and coarsening. The theoretical foundations of dendritic solidification lie in the classical Stefan problem, a moving boundary problem, that describes the evolution of a solid-liquid phase front~\citep{rubinstein1971stefan}. The two-phase interface motion is obtained by solving the heat equation in each phase, coupled with an evolving interface boundary condition (Stefan condition) that explicitly sets the velocity of the moving interface. While the evolution of dendrites in pure melts with no dissolved solutes is well represented by the Stefan problem, for melts with dissolved solutes (like alloys), mass diffusion in the constituent phases should also be accounted for in the governing equations. 
\par
Many analytical and numerical approaches exist for treating the Stefan problem, albeit under various simplifying assumptions on the problem geometry, interface geometry and boundary conditions. In the numerical domain, front tracking methods~\citep{meyer1978numerical, marshall1986front} are useful in solving straight and curved interface evolution in one-dimensional problems. Other popular methods that explicitly track moving interface are the Landau transformation and the finite element mesh-based moving node techniques. These methods are more suited when the movement of the interface is not far from the initial position~\citep{dantzig2016solidification}, but bookkeeping of the interface movement can be an arduous task, especially for anisotropic interphase growth conditions. The level set method is another popular numerical technique that is extensively used to solve moving interface problems, including dendrite solidification using the Stefan problem~\citep{chen1997simple}. 
\par
The phase-field method has gained popularity in the last two decades as an alternate diffuse interface numerical technique to model solidification processes. The phase-field model modifies the equations represented by the Stefan problem by introducing an order parameter to distinguish between the solid and liquid phases. Further, the previously sharp interface is represented as a diffuse interface using this order parameter without sacrificing much of the accuracy. The earliest of the simplified isotropic phase-field model has been proposed by Fix~\cite{fix1983phase}, and Collins and Levine~\cite{collins1985diffuse}. 
Caginalp~\cite{caginalp1986analysis, caginalp1989stefan} and  Kobayashi~\cite{kobayashi1993modeling} introduced basic anisotropy into the phase-field model and numerically predicted dendritic patterns seen in the solidification of the pure melt. They also modeled the formation of the side branching by introducing noise into the phase-field models.  
\par
The potential of the so far surveyed phase-field models was limited due to constraints on the lattice size, undercooling, interface width, capillary length, and interface kinetics. A modified thin interface phase-field model proposed by Karma and Rappel~\cite{karma1996phase} overcame these limitations. The dendrite shape and tip velocities of the pure melt obtained from the thin interface phase-field model were in excellent agreement with the steady-state results reported in theoretical studies. Karma and Rappel~\cite{karma1999phase} studied dendrite side-branching by introducing a microscopic thermal noise, and the side branch characteristics were in line with the linear WKB theory based stability analysis.  Plapp and Karma~\citep{plapp2000multiscale} proposed a hybrid technique where the phase-field model and diffusion-based Monte-Carlo algorithm were used in conjunction to simulate dendritic solidification at the low undercooling in two and three-dimensional geometry. 
\par
The phase-field model for the solidification of a binary alloy was presented by Warren and Boettinger~\cite{warren1995prediction}. Realistic patterns showing the primary and the secondary arms of the dendrite were simulated using the interface thickness far different from the asymptotic limit of the sharp interface models. Loginova \textit{et al.}~\cite{loginova2001phase} was the first to include the non-isothermal effects in binary alloy solidification. However, the computations were time-consuming. An artificial solute trapping was observed in the simulations with higher undercooling. This artificial solute trapping was the result of diffuse interface thickness. Karma~\cite{karma2001phase} developed a modified thin interface phase-field model for alloy solidification. Compared to the previous models, this model was applicable for varied interface thicknesses, zero non-equilibrium effects at the interface, and nearly zero diffusivity in the solid.
\par
The thin interface phase-field model by Ramirez \textit{et al.}~\cite{ramirez2004phase} considered both heat and solute diffusion along with zero kinetics at the interface. The interface thickness used was an order of magnitude less than the radius of curvature of the interface.  Echebarria \textit{et al.}~\cite{echebarria2004quantitative} modeled numerous numerical test cases and studied the convergence of the thin interface phase-field model. They discussed the final form of anti-trapping current that is often used in many alloy solidification models~\citep{karma2001phase, ramirez2004phase, almgren1999second}. When the solid diffusivity of the alloys is non-negligible, the form of anti-trapping current previously used needs modification and this was rigorously derived  in the work of  Ohno and Matsuura~\cite{ohno2009quantitative}

\par
With the substantial review tracing Stefan problem and important theoretical development of the phase-field  model,  we review literature focusing on the numerical implementations of these models and error analysis.  Several important studies leveraged the finite-element and finite-difference method based numerical schemes and their application to the phase-field model of solidification.  These studies focused on the accuracy and stability of proposed numerical schemes.  In this regard,  the work of  Feng and Prohl~\citep{feng2004analysis} is important.  They made use of fully discrete finite element methods and obtained optimal error bounds in relation to the interface thickness for a phase-field method of solidification.  Using the error estimates they showed convergence of finite element schemes to the solution to the phase field models in the limit of the sharp interface.  Gonzalez-Ferreiro \textit{et al.}~\citep{gonzalez2014thermodynamically}  presented a finite element discretization in space and mid-point discretization in time of a phase-field model that is consistent with both the thermodynamics laws.  The thermodynamically consistent numerical model was developed to have better dendrite resolution and higher order accuracy in time.  A first-order accurate in time and energy stable numerical method was presented by Chen and Yang~\cite{chen2019efficient}. The proposed method was applied to coupled Allen-Cahn, heat diffusion, and modified Navier-Stokes equations.  They resort to techniques that decoupled these three equations with the use of implicit-explicit schemes in their numerical implementation.  
\par
Better spatial and temporal accuracy and smaller error estimates were reported in the work of Kessler and Scheid~\citep{kessler2002priori}. They applied a finite element method to the phase-field model of binary alloy solidification where the error convergence numerical tests were done on a physical example of Ni-Cu alloy solidification with a simpler and non-branched solidification structure.   The adaptive meshing technique was used in many numerical studies of solidification using a phase-field model.  Hu \textit{et al.}~\citep{hu2009multi} presented a multi-mesh adaptive finite-element based numerical scheme to solve the phase-field method for the pure-melt solidification problem. The accuracy of their method was tested with its ability to predict dendrite tip velocities with a range of undercooling conditions.  Work of Rosam \textit{et al.}~\citep{rosam2007fully} using the mesh and time adaptive fully implicit numerical scheme utilizing finite-difference method applied to binary alloy solidification is very relevant.  Such highly space-time adaptive techniques saved computational time.  Their focus was on the implicit and explicit schemes error estimates using only derived variables such as dendrite tip position,  radius, and velocity.  However, error analysis reported in the literature did not explore the effect of basis continuity (i.e, $C^n$-continuous basis) on capturing dendrite kinetics and morphologies, and the convergence studied reported did not explicitly study the error in the primal fields, i.e., temperature and phase-field order parameter. As part of this manuscript, these two aspects of error analysis will also be addressed.

A review by Tourret \textit{et al.}~\citep{tourret2022phase} highlights several studies that applied the phase-field methods to model solidification and obtained good comparisons with the experimental results.  Present challenges and extensions of phase-field models were also discussed in their review. Wang \textit{et al.}~\citep{wang2012phase} discovered the relation between lower and upper limit primary dendrite arm spacing to inter-dendritic solute distribution and inter-dendritic undercooling by solving phase-field models using finite-element method. Fallah \textit{et al.}~\citep{fallah2012phase} work showed the suitability of phase-field models coupled with heat transfer models to reproduce experimentally known complex dendrite morphology and accurate dendrite size of Ti-Nb alloys thereby demonstrating the suitability of the numerical methods to laser deposition and industrial scale casting process.   More recently, microstructure evolution processes involving competitive grain growth of columnar dendrites and the grain growth along the converging grain boundary were modeled using the phase-field model by Tourret and Karma~\cite{tourret2015growth}, and Takaki \textit{et al.}~\cite{takaki2014two}.  In manufacturing processes like welding and molding, phase-field models were used to study solidification cracking susceptibility in Al-Mg alloys by Geng \textit{et al.}~\citep{geng2018effects}, and dendrite morphology in the melt pool of Al-Cu alloys by Farzadi \textit{et al.}~\citep{farzadi2008phase}. Integrated phase-field model and finite element methods are also used to study microstructure evolution  especially formation of laves phase in additive manufacturing of IN718~\cite{wang2019investigation}.  Phase-field models finds extensive applications in modeling rapid solidification conditions prevalent in the metal additive manufacturing~\cite{rolchigo2017modeling, ghosh2017primary, gong2015phase, sahoo2016phase, keller2017application}. 

A lot of existing studies on solidification and related phase-field models have focused only on either pure metals or alloys.  However, in this work, we present a unified treatment of both pure metals and alloys. We discuss classical Stefan problems relevant to both these types of solidification, and present their phase-field formulations and numerical implementations. Further, we present novel numerical error estimates of dendritic tip velocity, and the convergence of error for the primal fields of temperature and order parameter with respect to the numerical discretization. Lastly, using this numerical framework, various types of physically relevant dendritic solidification patterns like single equiaxed, multi-equiaxed, single columnar and multi-columnar dendrites are modeled in two-dimensional and three-dimensional computational domains. 
\par  

Following this literature review and summary of the advancements in numerical modeling of solidification and dendritic growth, we now present an overview of this manuscript.  In Section~\ref{sub:NumericalModel} we present a detailed discussion of the numerical models used in this work to model dendritic growth. Essentially, we look at the appropriate Stefan problem formulations and their phase-field representations for modeling the solidification problem of pure melts and binary alloys.  Then, in Section~\ref{sub:numericalImplementation}, we discuss the numerical framework and its computational implementation. Further, we present a numerical error analysis and the simulation results of various 2D and 3D dendritic solidification problems. Finally, we present the concluding remarks in Section~\ref{sub:DiscussionConclusion}. 
   

\section{Numerical models of dendritic growth}\label{sub:NumericalModel}
In this section, we discuss the Stefan problem for modeling solidification of a pure metal and a binary alloy. We elaborate on the complexities of the Stefan problem and the physics associated with solidification. The Stefan problem is then re-written as a phase-field formulation. The primary goal of using a phase-field approach is to model the interface growth kinetics, particularly the velocity of interface motion and the complex morphologies of dendrites that occur during solidification, without the need to explicitly track the dendritic interfaces. The interface growth can be influenced by the surface anisotropy, heat diffusion, mass diffusion, interface curvature, and the interface attachment kinetics. In Sections~\ref{subsub:DetailedStefanPureMetal}-\ref{subsub:PhasefieldPureMetal}, the classical Stefan problem and its phase-field representation relevant to solidification of pure melts is presented. Then, in Sections~\ref{subsub:StefanBinaryAlloy}-\ref{subsub:PhasefieldAlloy}, the extensions needed in the Stefan problem and its phase-field representation for modeling solidification of binary alloys are discussed. 

In general, a pure melt is a single-component material (pure solvent) without any solutes. Addition of one or more solute components results in an alloy. In this work, we specifically consider a binary alloy (two-component alloy with a solvent and one solute), as a representative of multi-component alloys. And as will be shown, during the solidification of a binary alloy, we model the diffusion of this solute concentration in the solid and liquid phases.

\subsection{Stefan problem for modeling solidification of a pure melt} {\label{subsub:DetailedStefanPureMetal}}
The Stefan problem describing the solidification of an undercooled pure metal is presented. The set of Equations~\ref{eq:Stefan_T}-\ref{eq:Stefan_Interface} mathematically model the solidification of an undercooled pure melt.  Equations~\ref{eq:NDStefan_T}-\ref{eq:NDStefan_Interface} represent the same governing equations in a non-dimensional form.  Equation~\ref{eq:Stefan_T} models heat conduction in the bulk liquid and the bulk solid regions of the pure metal. At the solid-liquid interface, the balance of the heat flux from either of the bulk regions is balanced by the freezing of the melt, and thus additional solid bulk phase is formed. In other words, the solid-liquid interface moves. The energy balance at the interface is given by the Equation~\ref{eq:Stefan_Condition_T}. The temperature of the interface is not fixed and can be affected by multiple factors such as the undercooling effect on the interface due to its curvature (Gibbs-Thomson effect) and the interface attachment kinetics. These are captured in the equation~\ref{eq:Stefan_Interface}.  
\begin{subequations}\label{eq:Stefan_PureMetal}
\begin{equation} \label{eq:Stefan_T}
 \frac{\partial T(\boldsymbol{x},t)}{\partial t}=\frac{k}{c_p}\nabla^2 T,  \quad  \boldsymbol{x} \in \Omega^{s},  \Omega^{l}
\end{equation}
\begin{equation} \label{eq:Stefan_Condition_T}
\boldsymbol{\nu}_{n}L=k \big(\partial_{n}T(\boldsymbol{x})\lvert^{+} -\; \partial_{n}T(\boldsymbol{x})\lvert^{-} \big), \quad \boldsymbol{x} \in \Gamma
\end{equation}
\begin{equation}\label{eq:Stefan_Interface}
T_m-T_i=\Delta T_{\Gamma_g} +\Delta T_{\mu_k}=\Gamma_g\boldsymbol{\kappa} + \frac{\boldsymbol{\nu}_{n}}{\mu_k} , \quad \boldsymbol{x} \in \Gamma
\end{equation}
\end{subequations}
Here, $k$ and $c_p$ are thermal conductivity and specific heat capacity for a pure solid or liquid metal. $T_m$ $\boldsymbol{\nu}_n$, $L$ is the melting temperature, interface velocity, and latent heat of the metal. $\partial_{n}T\lvert^{+}$ and $\partial_{n}T\lvert^{-}$ is the temperature gradient normal to the interface in the liquid and the solid phase respectively. $\Gamma_g$, $\boldsymbol{\kappa}$, and $\mu_k$ are the Gibbs-Thomson coefficient, the curvature of the interface, and interface-attachment coefficient, respectively.  We re-write these equations in their non-dimensional form, using the scaled temperature $u= \frac{c_p(T-T_m)}{L}$ \citep{dantzig2016solidification}. 
\begin{subequations}\label{eq:Stefan_PureMetalND}
\begin{equation} \label{eq:NDStefan_T}
 \frac{\partial u (\boldsymbol{x},t) }{\partial t}=D\nabla^2 u ,  \quad  \boldsymbol{x} \in \Omega^{s},  \Omega^{l}
\end{equation} 
\begin{equation} \label{eq:NDStefan_Condition_T} 
\boldsymbol{\nu}_{n}=D \big(\partial_{n}u(\boldsymbol{x})\lvert^{+} -\; \partial_{n}u(\boldsymbol{x})\lvert^{-} \big), \quad \boldsymbol{x} \in \Gamma
\end{equation}
\begin{equation}\label{eq:NDStefan_Interface}
u^{*}=-d(\mathbf{n}) \boldsymbol{\kappa} - \beta(\mathbf{n})\boldsymbol{\nu}_{n}, \quad \boldsymbol{x} \in \Gamma
\end{equation}
\end{subequations}
$D=\frac{k}{c_p}$ is the non-dimensional thermal diffusivity in the solid and liquid phases. 
$d(\mathbf{n})=\gamma(\mathbf{n})T_{m}c_{p}/L^{2} $ is the capillary length, $\gamma(\mathbf{n})$ is the surface tension, $\mathbf{n}$ is the unit vector denoting normal to the interface,   $\beta (\mathbf{n})$ is the kinetic coefficient, and $\partial_{n}u\lvert^{+}$ and $\partial_{n}u\lvert^{-}$ are the non-dimensional derivatives normal to the interface in the solid and liquid region, respectively \citep{slavov2003phase}. The driving force for this solidification is the initial undercooling, i.e., the initial temperature of the liquid melt below its melting temperature. The non-dimensional undercooling is given by $\Delta=-c_p(\frac{Tm-T_\infty}{L})$. 

The Stefan problem described using Equation~\ref{eq:Stefan_PureMetal} does not have a known analytical solution, and its numerical implementation, as is the case with free boundary problems, is challenging. Numerical methods like time-dependent boundary integral formulations \citep{karma1996phase}, a variational algorithm with zero or non-zero interface kinetics, and methods that require the book-keeping of the solidifying front and low grid anisotropy conditions have been used in the past. Modeling correct dendritic growth involves getting the correct operating dendritic tip conditions, mainly the tip radius and the tip velocity. Variationally derived phase-field models, like the one described in Section~\ref{subsub:PhasefieldPureMetal}, are now the method of choice to model complex dendritic solidification problems.  

\subsection{Phase-field model describing the solidification of a pure melt}{\label{subsub:PhasefieldPureMetal}}
In this section, we describe the phase-field method developed by Plapp and Karma~\cite{plapp2000multiscale} to model solidification in a pure metal. The governing equations for the model are derived using variational principles. A phenomenological, \textcolor{black}{isothermal free energy expression, $\Pi \left[\phi, u\right]$, in terms of a phase-field order parameter ($\phi$) and temperature ($u$), is given by Equation~\ref{eq:FunctionalForm}. }.This functional form defines the thermodynamic state of a system, and has two terms - first is $f(\phi,u)$, the bulk energy term, and second is an interface term given by $\frac{1}{2} \lambda^2(\boldsymbol{n})| \bm{\nabla} \phi|^2$. The magnitude of the interface term is controlled by the interface parameter, $\lambda$, and the interface term is positive in the diffuse solid-liquid interface due to a non-zero gradient of the order parameter; elsewhere, the interface term is zero. 
\begin{subequations}\label{eq:FunctionalForm} 
\begin{equation} \label{eq:FunctionalPureMetal}
  \Pi \left[\phi, u\right] = \int_{\Omega}   \left[ f(\phi,u) + \frac{1}{2} \lambda^2(\mathbf{n})| \bm{\nabla} \phi|^2 \right]   ~dV 
\end{equation}
\begin{equation}\label{eq:doubleWell}
f(\phi,u) = -\frac{1}{2}\phi^2 + \frac{1}{4}\phi^4 + \xi u \phi \left(1-\frac{2}{3} \phi^2+\frac{1}{5}\phi^4 \right)
\end{equation}
\end{subequations}

\noindent The bulk energy term $f(\phi,u)$ in Equation~\ref{eq:doubleWell} is given by the double well potential that has a local minima at $\phi=1$ and $\phi=-1$.  $\Omega = \Omega_s \cup \Omega_l $ represents total volume encompassing the bulk solid,  the bulk liquid and the interface region. The value of $u$ then tilts the equilibrium value of $f(\phi,u)$. The form of $f(\phi,u)=g(\phi)+\xi u h(\phi) $ is given below. $\xi$ is the coupling parameter in the double well function. The governing equations for the phase-field model of solidification are obtained by minimizing the above functional with respect to the primal fields, $\phi$.   \textcolor{black}{
The governing equations for the temperature field, $u$, are obtained by taking the variational derivative of the Lyapunov functional: $\Pi_L \left[\phi, u\right] = \int_{\Omega}   \left[ g(\phi) +\xi u^2 + \frac{1}{2} \lambda^2(\mathbf{n})| \bm{\nabla} \phi|^2 \right]   dV $, whose isothermal form is given by Equation~\ref{eq:FunctionalPureMetal}.   A detailed discussion on this functional and its isothermal form are beyond the scope of this manuscript, and interested readers are referred to the relevant literature ~\cite{plapp2000multiscale, boettinger2002phase, karma1998quantitative}. The non-dimensional enthalpy, $U$, accounts for the change in temperature and the latent heat, $U=u-\frac{1}{2}\phi$. In terms of the enthalpy, the governing equations for the first-order kinetics are given by:}
\begin{subequations}\label{eq:PhaseFieldUandPhiShort}
\begin{equation} \label{eq:VD_T}
\textcolor{black}{
\frac{\partial U(\boldsymbol{x},t)}{\partial t}= \boldsymbol{\nabla}. \Big ( D \boldsymbol{\nabla} \frac{\delta \Pi_L}{\delta u}\Big )} ,  \quad  \boldsymbol{x} \in \Omega
\end{equation}
\begin{equation} \label{eq:VD_PPhi}
\textcolor{black}{
\tau(\bm{n})\frac{\partial \phi (\boldsymbol{x},t) }{\partial t}= -\frac{\delta \Pi}{\delta \phi} } ,  \quad  \boldsymbol{x} \in \Omega  
\end{equation}
\end{subequations} 
\textcolor{black}{Rewriting these equations in terms of the primal fields, $u$ and $\phi$, yields~\cite{plapp2000multiscale},}
\begin{subequations}\label{eq:PhaseFieldUandPhi}
\begin{equation} \label{eq:EquationU}
 \frac{\partial u }{\partial t}=D\nabla^2 u + \frac{1}{2}\frac{\partial\phi}{\partial t}, \quad  \boldsymbol{x}  \in \Omega   
\end{equation}
\begin{align} \label{eq:EquationPHI}
\tau(\boldsymbol{n}) \frac{\partial \phi}{\partial t} = -\frac{\partial f}{\partial \phi} + \bm{\nabla} \cdot \big(\lambda^2(\boldsymbol{n}) \bm{\nabla} \phi \big)  +\frac{\partial}{\partial \bm{x}} \left[ | \bm{\nabla} \phi|^2 \lambda(\boldsymbol{n}) \frac{\partial \lambda(\boldsymbol{n})}{\partial \left( \frac{\partial \phi}{\partial \bm{x}} \right)} \right] + \frac{\partial}{\partial \bm{y}} \left[ | \bm{\nabla} \phi|^2 \lambda(\boldsymbol{n}) \frac{\partial \lambda(\boldsymbol{n})}{\partial \left( \frac{\partial \phi}{\partial \bm{y}} \right)} \right],  \boldsymbol{x}  \in \Omega  
\end{align}
\end{subequations} 
 
\noindent where $\Omega$ denotes the problem domain. Anisotropy in the interface energy representation is considered. Thus, $\lambda(\boldsymbol{n})=\lambda_{0}a_{s}(\boldsymbol{n})$ and \textcolor{black}{$\tau(\boldsymbol{n})=\tau_{0}a_{s}^2(\boldsymbol{n})$, where $\tau_{0}$ is the characteristic time scale of evolution of the order parameter.} The anisotropy parameter is defined as  $a_{s}=(1-3\epsilon_{4})\Big [ 1+ \frac{4\epsilon_{4}}{(1-3\epsilon_{4})}\frac{(\partial_{x}\phi)^4+(\partial_{y}\phi)^4}{|\phi|^4} \Big ]=1+\epsilon_4cos(m\theta)$. The equivalence is easily established by taking $\tan(\theta) = \frac{\partial \phi}{\partial y} / \frac{\partial \phi}{\partial x}$. 
$\epsilon_4$ is the strength of anisotropy and m=4 corresponds to an anisotropy with four-fold symmetry. The interface between the solid and liquid phases is diffuse in the phase-field representation of a Stefan problem of solidification, but in the asymptotic limit of  $\lambda(\boldsymbol{n}) \to 0$, the Stefan problem is recovered~\citep{caginalp1986analysis, caginalp1989stefan}.

\subsubsection{Weak formulation}
We now pose the above governing equations in their weak (integral) form. This formulation is used to solve these equations within a standard finite element method framework.\\ \\ Find the primal fields $\{ u, \phi \}$, where,
\begin{align*}
u &\in \mathscr{S}_{u},  \quad \mathscr{S}_{u} = \{ u \in \text{H}^1(\Omega)~\vert  ~u  = ~u' ~\forall ~\textbf{X} \in \Gamma^{u} \}, \\
\phi&\in \mathscr{S}_{\phi},  \quad \mathscr{S}_{\phi} = \{ \phi \in \text{H}^1(\Omega)~\vert  ~\phi = ~\phi' ~\forall ~\textbf{X} \in \Gamma^{\phi} \}
\end{align*}
such that,   
\begin{align*}
\forall ~w_{u} &\in \mathscr{V}_{u},  \quad \mathscr{V}_{u} = \{ w_{u} \in \text{H}^1(\Omega)~\vert  ~w_{u}  = 0 ~\forall ~\textbf{X} \in \Gamma^{u} \}, \\
\forall ~w_{\phi} &\in \mathscr{V}_{\phi},  \quad \mathscr{V}_{\phi} = \{ w_{\phi} \in \text{H}^1(\Omega)~\vert  ~w_{\phi}  = 0 ~\forall ~\textbf{X} \in \Gamma^{\phi} \}
\end{align*}
we have, 

\begin{subequations}\label{eq:WeakFormPhaseFieldUPhi}
\begin{equation}\label{eq:WeakFormU}
\int_{\Omega}w_u \Big( \frac{\partial u }{\partial t} - \frac{1}{2}\frac{\partial\phi}{\partial t}  \Big)  \,dV + \int_{\Omega}D \big (\bm{\nabla} w_u . \bm{\nabla} u \big)  \, dV  = 0
\end{equation}
\begin{align}\label{eq:WeakFormPhi}
\int_{\Omega} w_{\phi} \Big (  \tau \frac{\partial \phi}{\partial t}  + \frac{\partial f}{\partial \phi} \Big )\,dV + \int_{\Omega} \lambda^2 \bm{\nabla} w_{\phi} .  \bm{\nabla} \phi \,dV +  \int_{\Omega} | \bm{\nabla} \phi|^2 \lambda \frac{ \partial w_{\phi} }{\partial \bm{x}} \frac{\partial \lambda}{\partial \left( \frac{\partial \phi}{\partial \bm{x}} \right)}\,dV \\ \nonumber +  \int_{\Omega} |\bm{\nabla} \phi|^2 \lambda \frac{ \partial w_{\phi} }{\partial \bm{y}} \frac{\partial \lambda}{\partial \left( \frac{\partial \phi}{\partial \bm{y}} \right)}\,dV
= 0
\end{align}
\end{subequations} 
here, $\Omega = \Omega_s \cup \Omega_l $, where $\Omega_s$ is the solid phase and $\Omega_l$ is the liquid phase, respectively. The variations, $\omega_{u}$  and $\omega_{\phi}$, belong to $\text{H}^1(\Omega)$ - the Sobolev space of functions that are square-integrable and have a square-integrable derivatives.
 
\subsection{Stefan problem for modeling solidification of a binary alloy} \label{subsub:StefanBinaryAlloy}
\begin{figure}[ht!]
  \centering
  \subfloat[\label{fig:Phase_Diagram}]{
                \includegraphics[width=0.4\textwidth]{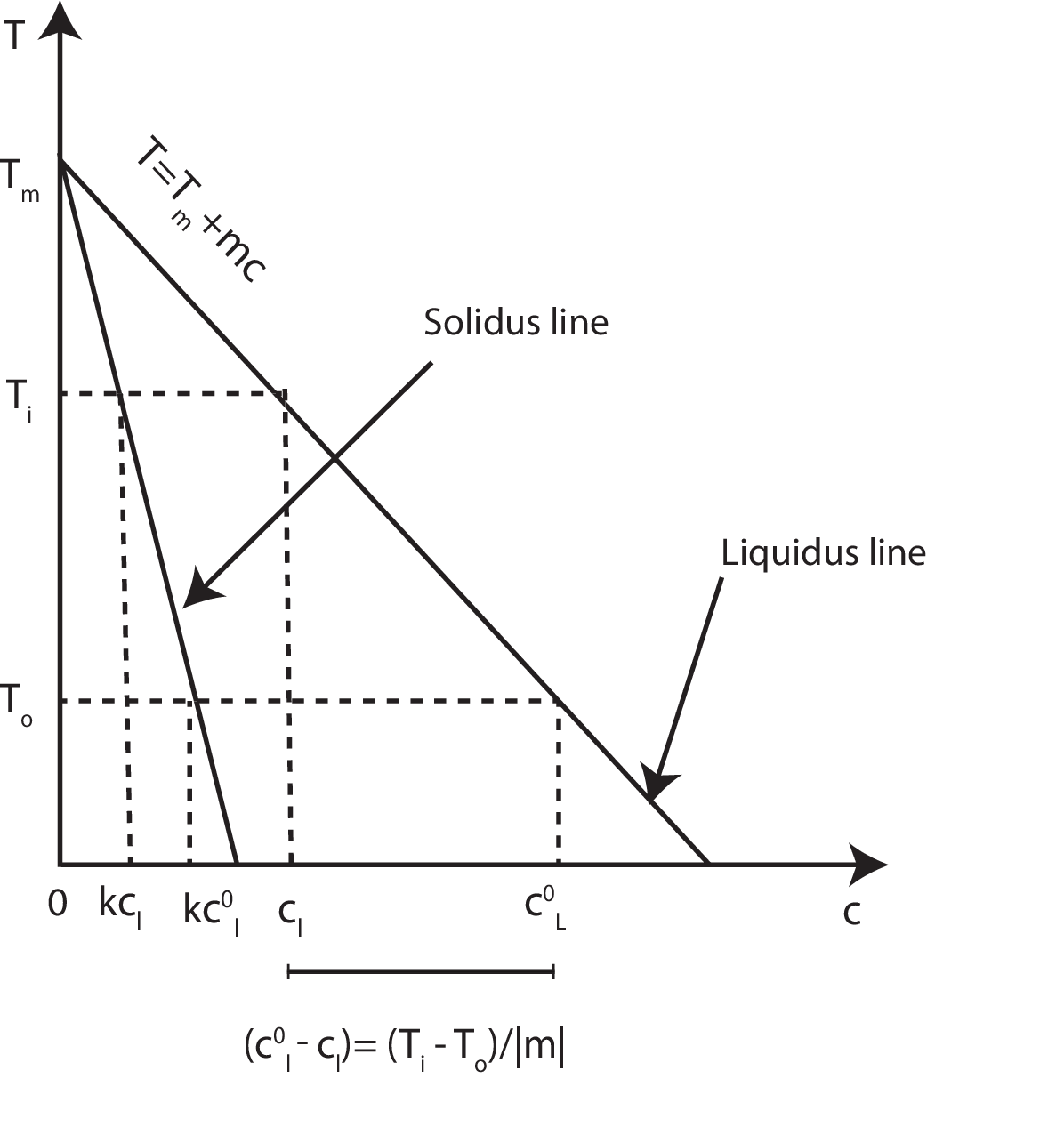} 

  }
   \subfloat[\label{fig:Directional}]{
                \includegraphics[width=0.4\textwidth]{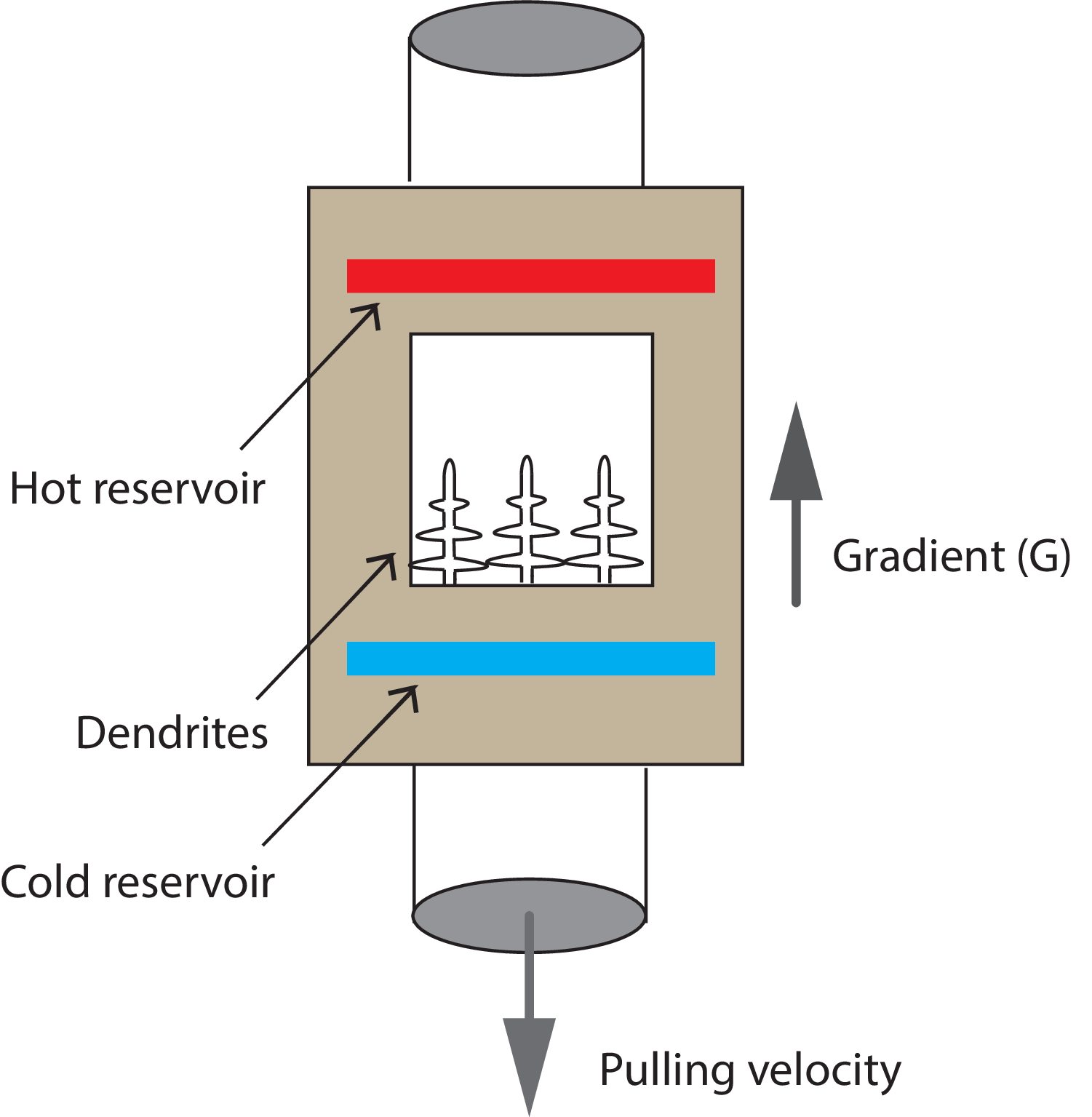}

  }
  \caption{Directional solidification in a binary alloy. Shown are (\ref{fig:Phase_Diagram}) a simplified phase diagram of a binary alloy, and (\ref{fig:Directional}) a schematic of a directional solidification process. In the case of directional solidification, typically, the melt is forced to solidify at a constant pulling velocity under a fixed temperature gradient. The dendritic structures in this case evolve at the solid-liquid interface, via a Mullins-Sekerka instability, when the pulling velocity is above a critical value.}
  \label{fig:Phase_Directional}
\end{figure}

The Stefan problem describing the solidification of a binary alloy is described in this section. As mentioned earlier, an alloy is a multi-component system. A binary alloy is a specific case of a multi-component alloy with only one solute component. A binary alloy consists of a solute that has a composition of $c$, and a solvent that has a composition of $(1-c)$. During the solidification of a binary alloy, we model the diffusion of the solute concentration across the solid and the liquid phase. Mass diffusion is governed by Equation~\ref{eq:Stefan_C}.
Near the vicinity of the solid-liquid interface, the composition of the solute in the solid and the liquid phase is dictated by the binary alloy phase diagram given in the Figure~\ref{fig:Phase_Diagram}. The difference in the concentration of the solute across the solid and liquid phase is balanced by the mass diffusion flux. This is captured in Equation~\ref{eq:Stefan_Condition_C}. The temperature of the interface is not fixed and can be affected by multiple factors.  The undercooling on the interface can be due to its curvature (Gibbs-Thomson effect),  the interface attachment kinetics, cooling rate ($\dot{T}$), and the applied temperature gradient. This is mathematically represented in Equation~\ref{eq:Stefan_Interface_C}.  Using the phase-diagram, we can write this condition in-terms of composition, $c_l$ at the interface. We neglect the contribution from the attachment kinetics. $c_l=c^0_l - \frac{\Gamma_g \boldsymbol{\kappa} + Gy + \dot{T}t}{|m|} $ where $c$ is the solute composition in either the solid or the liquid phase. $D$ is the mass diffusivity of alloy in the liquid phase, $c_l$ and $c^0_l$ are the solute concentration at the interface for the temperature $T_{i}$ and $T_{0}$ respectively. $\dot{T}=G\bm{\nu}_{p}$ is the cooling rate. We scale the solute concentration into a variable $u_c=\frac{c-c^0_l}{c^0_l(1-k)}$ which represents supersaturation. k is the partition coefficient which represents the ratio of solid composition to the liquid composition.  
\begin{subequations}\label{eq:Stefan_Alloy}
\begin{equation} \label{eq:Stefan_C}
 \frac{\partial c (\boldsymbol{x},t) }{\partial t}=D\nabla^2 c (\boldsymbol{x},t),  \quad  \boldsymbol{x} \in \Omega^{s},  \Omega^{l}
\end{equation} 
\begin{equation} \label{eq:Stefan_Condition_C}
 c_{l}(1-k) \boldsymbol{\nu}_{n}=-D\partial_{n}c (\boldsymbol{x},t) \lvert^{+}, \quad  \boldsymbol{x} \in \Gamma
\end{equation}
\begin{equation}\label{eq:Stefan_Interface_C}
T_m-T_i=\Gamma \boldsymbol{\kappa} + \frac{ \boldsymbol{\nu}_n}{\mu_k}+ Gy + \dot{T}t, \quad  \boldsymbol{x} \in \Gamma
\end{equation}
\end{subequations}
\noindent here, $\Gamma$ is the solid-liquid interface that separates the solid phase, $\Omega^{s}$, and the liquid phase, $\Omega^{l}$. \\

\noindent Rewriting the Equations~\ref{eq:Stefan_C}, \ref{eq:Stefan_Condition_C} and \ref{eq:Stefan_Interface_C} in-terms of the supersaturation, $u_c$, we get, 
\begin{subequations}\label{eq:NDStefan_Alloy}
\begin{equation} \label{eq:NDSTE_U}
 \frac{\partial u_{c} (\boldsymbol{x},t) }{\partial t}=D\nabla^2 u_{c} (\boldsymbol{x},t), \quad  \boldsymbol{x} \in \Omega^{s},  \Omega^{l}
\end{equation}
\begin{equation} \label{eq:NDStefan_Condition_C}
 (1+(1-k)u^{*}_{c}) \boldsymbol{\nu}_{n} =-D\partial_{n}u_{c}(\boldsymbol{x},t)\lvert^{+},  \quad  \boldsymbol{x} \in \Gamma
\end{equation}
\begin{equation} \label{eq:NDStefan_Interface_C}
u^{*}_{c}=-d_{0} \boldsymbol {\kappa} + \dot{\theta}t-\gamma z,  \quad  \boldsymbol{x} \in \Gamma
\end{equation}
\end{subequations}

\noindent Here $d_{0}=\frac{\Gamma_g}{\Delta T_{0}}$ is the chemical capillary length , $\dot{\theta}=\frac{-\dot{T}}{\Delta T_{0}}$ is the  non-dimensional cooling rate and $\gamma=\frac{G}{\Delta T_{0}}$ is non-dimensional thermal gradient and $\Delta T_{0}=|m|(1-k)c^{0}_{l}$ \citep{neumann2017phase}. 

\noindent As can be expected, an analytical solution for the moving boundary problem in Equations~\ref{eq:NDStefan_Alloy} is not known. In the next subsection, we discuss a phase-field representation of this problem that is extensively used to model dendritic solidification of binary alloys.


\subsection{Phase-field model describing the solidification of a binary alloy}{\label{subsub:PhasefieldAlloy}}
In this section, we describe the phase-field model for the solidification of a binary alloy. Following the variational procedure adopted in the case of a pure melt, we write the free energy functional form, $\Pi[\phi,c, T]$, for a binary alloy. The free energy in this case has the additional dependence on the solute composition, $c$, along with the order paramter, $\phi$, and the temperature, $T$. This functional form, given by Equation~\ref{eq:FunctionalAlloy}, takes into account bulk free energy of the system $f(\phi,T_m)$, internal energy and entropy $f_{AB}(\phi,c,T)$, and the interface term  $\frac{1}{2} \sigma|\bm{\nabla} \phi|^2$, where  $\sigma \propto \lambda^2$.  
\begin{equation} \label{eq:FunctionalAlloy}
  \Pi[\phi,c,T] = \int_{\Omega}   \left[ f(\phi,T_{m}) + f_{AB}(\phi,c,T) + \frac{1}{2} \sigma|\bm{\nabla} \phi|^2 \right]   ~dV 
\end{equation}
\noindent here $\Omega = \Omega_s \cup \Gamma \cup \Omega_l $ represents the total volume encompassing the bulk solid, the bulk liquid and the interface region.  
The governing equations for the model are obtained by writing the diffusion equations relevant to a conserved quantity - the solute composition, and a non-conserved quantity - the order-parameter. 
\begin{subequations}\label{eq:PhaseFieldUandPhiAlloy}
\begin{equation} \label{eq:VD_C}
\frac{\partial c (\boldsymbol{x},t)}{\partial t}= \bm{\nabla}. \Big ( M(\phi,c) \bm{\nabla} \frac{\delta \Pi}{\delta c}\Big ), \quad  \boldsymbol{x} \in \Omega
\end{equation}
\begin{equation} \label{eq:VD_phi}
\frac{\partial \phi (\boldsymbol{x},t)}{\partial t}= -K_{\phi}\frac{\delta \Pi}{\delta \phi}, \quad  \boldsymbol{x} \in \Omega
\end{equation}
\noindent where $M(\phi,c) $ is the solute mobility, and $K_{\phi}$ is the order parameter mobility. It is convenient to scale the composition, c and write it as $u_c=\frac{1}{1-k}\big ( \frac{2c}{c^0_l(1-\phi + k (1+\phi))}-1 \big )$. The parameter $u_c$ can be understood as a concentration undercooling parameter analogous to $u$ in Section~\ref{subsub:PhasefieldPureMetal}. This change of variable enables us to draw similarities with the earlier form of the phase-field model described in Equations~\ref{eq:PhaseFieldUandPhi}.   \textcolor{black}{At the onset of solidification, the initial solute concentration in the liquid region is $c=c^0_l$ and the solute concentration in the solid region is $c=kc^0_l$, as per the phase-diagram.  However,  when expressed in term of $u_c$, it is instructional to note that the initial composition undercooling,  $u_{c}=0$, both in the solid region (${c=kc^0_l, \phi=1}$) and the liquid region (${c=c^0_l, \phi=-1}$).  Thus, it is preferable to work with $u_{c}$ as the primal field, as it is a continuous variable across the solid and liquid region unlike the solute composition, $c$, which is discontinuous, as can be easily seen from the initial conditions}.   We present the final form of the phase-field model in Equations~\ref{eqn:PF_Alloy_UC}-\ref{eqn:PF_Alloy_Phi}, without the accompanying mathematical derivation. For interested readers, a detailed derivation can be found in \cite{echebarria2004quantitative}.   
\begin{align} \label{eqn:PF_Alloy_UC}
    \Big(\frac{1+k}{2}-\frac{1-k}{2}\phi \Big)\frac{\partial u_c}{\partial t} = \bm{\nabla} \cdot \Big(\Tilde{D} \big (\frac{1-\phi}{2}\big) \bm{\nabla} u_c - \boldsymbol{J}_{at} \Big) + \Big(1+(1-k)u_c \Big)\frac{1}{2}\frac{\partial \phi}{\partial t},  \quad  \boldsymbol{x} \in \Omega
\end{align}
\begin{align} \label{eqn:PF_Alloy_Phi}
\begin{split}
\Big(1-(1-k)\frac{y-\tilde {\bm{\nu}} t}{\tilde l_{T}}\Big)a^2_{s}(\mathbf{n}) \frac{\partial \phi}{\partial t} = \bm{\nabla} \cdot \left(a^2_{s}(\mathbf{n}) \bm{\nabla} \phi \right)+  \frac{\partial}{\partial \bm{x}} \Big( | \bm{\nabla} \phi|^2 a_{s}(\mathbf{n}) \frac{\partial a_{s}(\mathbf{n})}{\partial \left( \frac{\partial \phi}{\partial \bm{x}} \right)} \Big) \\ + \frac{\partial}{\partial \bm{y}} \Big( | \bm{\nabla} \phi|^2 a_{s}(\mathbf{n}) \frac{\partial a_{s}(\mathbf{n})}{\partial \left( \frac{\partial \phi}{\partial \bm{y}} \right)} \Big) 
 +\phi-\phi^{3}-\xi(1-\phi^{2})^{2}\big (u_{c}+\frac{y-\tilde{\bm{\nu}}t}{\tilde l_{T}} \big ) 
\end{split}, \quad  \boldsymbol{x} \in \Omega
\end{align}

\noindent The solute diffusivity in the solid region is small as compared to the solute diffusivity in the liquid region. Thus mass diffusion is neglected in the solid region. $\tilde D$ is the non-dimensional diffusivity in the liquid phase. In the directional solidification process, dendrite growth is influenced by the presence of a temperature gradient $G$ and pulling velocity $\tilde{\bm{\nu}}_p$. This can be realized by the the term $\frac{y-\tilde{ \bm{\nu}} t}{\tilde l_{T}}$ in the Equation~\ref{eqn:PF_Alloy_Phi}, where $\tilde l_{T}= \frac{|m|(1-k)c^0_{l}}{G}$ is called the thermal length.  A modified form of Equation~\ref{eqn:PF_Alloy_UC} was suggested by Ohno and Matsuura\cite{ohno2009quantitative} where the mass diffusion in the solid phase is not neglected, and $D_{s}$ is comparable to $D_{l}$. This results in the modification to the diffusivity part and anti-trapping flux current as seen in the Equation~\ref{eqn:Alloy_Phi_MODIFIED}. If we substitute k=1 and $D_s=D_l$ in this equation, Equation \ref{eqn:PF_Alloy_UC} is recovered. 
\begin{equation}\label{eqn:Alloy_Phi_MODIFIED}
    \begin{split}
    \Big(\frac{1+k}{2}-\frac{1-k}{2}\phi \Big)\frac{\partial u_c}{\partial t} = \bm{\nabla} \cdot \Bigg(\Tilde{D} \big (\frac{1-\phi}{2}+k\frac{1+\phi}{2}\frac{D_{s}}{D_{l}}\big)\bm{\nabla u_c} - (1-k\frac{D_{s}}{D_{l}})\boldsymbol{J}_{at} \Bigg)\\
   + \Big(1+(1-k)u_c \Big)\frac{1}{2}\frac{\partial \phi}{\partial t}
    \end{split}
\end{equation}
\end{subequations}

\subsubsection{Weak formulation}

We now pose the above governing equations in their weak (integral) form. This formulation is used to solve these equations within a standard finite element framework. \\ \\ 
Find the primal fields $\{ u_c, \phi \}$, where,
\begin{align*}
u_c &\in \mathscr{S}_{u_c},  \quad \mathscr{S}_{u_c} = \{ u_c \in \text{H}^1(\Omega) ~\vert  ~u_c  = ~u_c' ~\forall ~\textbf{X} \in \Gamma^{u_c} \}, \\
\phi&\in \mathscr{S}_{\phi},  \quad \mathscr{S}_{\phi} = \{ \phi \in \text{H}^1(\Omega) ~\vert  ~\phi = ~\phi' ~\forall ~\textbf{X} \in \Gamma^{\phi} \}
\end{align*}
such that,   
\begin{align*}
\forall ~w_{u_c} &\in \mathscr{V}_{u_c},  \quad \mathscr{V}_{u_c} = \{ w_{u_c} \in \text{H}^1(\Omega)~\vert  ~w_{u_c}  = 0 ~\forall ~\textbf{X} \in \Gamma^{u_c} \}, \\
\forall ~w_{\phi} &\in \mathscr{V}_{\phi},  \quad \mathscr{V}_{\phi} = \{ w_{\phi} \in \text{H}^1(\Omega)~\vert  ~w_{\phi}  = 0 ~\forall ~\textbf{X} \in \Gamma^{\phi} \}
\end{align*}
we have,
\begin{subequations}\label{eq:WeakFormPhaseFieldUPhi_Alloy}
\begin{align}\label{eq:WeakFormU_Alloy}
\int_{\Omega}w_{u_{c}} \Bigg(  \Big(\frac{1+k}{2}-\frac{1-k}{2}\phi \Big)\frac{\partial u_c}{\partial t}  -  \Big(1+(1-k)u_c \Big)\frac{1}{2}\frac{\partial \phi}{\partial t}  \Bigg)  \, dV  \\ \nonumber 
+ \int_{\Omega} \Tilde{D} \big (\frac{1-\phi}{2}\big) \bm{\nabla} w_{u_{c}} .  \big( \bm{\nabla} u_c - \boldsymbol{J}_{at} \big)\,dV  = 0
\end{align}
\begin{align}\label{eq:WeakFormPhi_Alloy}
\int_{\Omega} w_{\phi} \Bigg( \Big(1-(1-k)\frac{y-\tilde{\bm{\nu}} t}{\tilde l_{T}}\Big)a^2_{s}(\mathbf{n}) \frac{\partial \phi}{\partial t}  - \Big ( \phi-\phi^{3}-\xi(1-\phi^{2})^{2}\big (u_{c}+\frac{y-\tilde{\bm{\nu}}t}{\tilde l_{T}} \big )  \Big)  \Bigg) \,dV  \\ \nonumber + \int_{\Omega} a_{s}^2 \bm{\nabla} w_{\phi} .  \bm{\nabla} \phi \,dV +  \int_{\Omega} | \bm{\nabla} \phi|^2 a_s \frac{ \partial w_{\phi} }{\partial \bm{x}} \frac{\partial a_s}{\partial \left( \frac{\partial \phi}{\partial \bm{x}} \right)}\,dV +  \int_{\Omega} | \bm{\nabla} \phi|^2 a_s \frac{ \partial w_{\phi} }{\partial \bm{y}} \frac{\partial a_s}{\partial \left( \frac{\partial \phi}{\partial \bm{y}} \right)}\,dV
= 0
\end{align}
here, $\Omega = \Omega_s \cup \Omega_l $, where $\Omega_s$ is the solid phase and $\Omega_l$ is the liquid phase, respectively. The variations, $\omega_{u_c}$  and $\omega_{\phi}$, belong to $\text{H}^1(\Omega)$ - the Sobolev space of functions that are square-integrable and have a square-integrable derivatives.

\end{subequations}

\section{Numerical implementation and results}\label{sub:numericalImplementation}

This section outlines the numerical implementation of the phase-field models described in Section~\ref{subsub:PhasefieldPureMetal} and \ref{subsub:PhasefieldAlloy}.  In Section~\ref{sub:ComputationalImp}, we cover the computational implementation. In Section~\ref{sub:InputParameters}, we summarize the input model parameters used in the simulation of dendritic growth in a pure metal and a binary alloy, and list the non-dimensional parameters used in these models. In Sections ~\ref{sub:ShapeSchematics}, using various dendritic morphologies modeled in this work, we identify important geometric features of dendrites that are tracked in the dendritic shape studies presented later. In Sections~\ref{sub:OptimalConvergence} we simulate the classical four fold symmetric dendrite shape occurring in undercooled pure melt and use it as a basis for convergence studies of dendritic morphology and dendritic tip velocity. Later, in Sections~\ref{sub:MultipleEquiaxedDendrite}-\ref{sub:3DSingleEquiaxed} , we model the evolution of complex 2D and 3D dendrite morphologies. These simulations cover various cases of solidification ranging from a single equiaxed dendrite in a pure liquid melt to a multi-columnar dendrites growing in a  binary alloy. Finally, evolution of equiaxed dendrite in 3D is demonstrated.

\subsection{Computational implementation} 
{\label{sub:ComputationalImp}}
The phase-field formulations presented in this work are solved using computational implementations of two numerical techniques: (1) A $C^0$-continuous basis code based on the standard Finite Element Method (FEM); the code base is an in-house, C++ programing language based, parallel code framework with adaptive meshing and adaptive time-stepping -  build on top of the deal.II open source Finite Element library~\cite{dealII93}, (2) A $C^1$-continuous basis code based on the Isogeometric Analysis (IGA) method; the code base is an in-house, C++ programing language based, parallel code framework -  build on top of the PetIGA open source IGA library~\citep{PetIGA}. As can be expected, the IGA code base can model $C^0$-continuous basis also, but does not currently support some capabilities like adaptive meshing that are needed for 3D dendritic simulations. Both codes support a variety of implicit and explicit time stepping schemes. Following the standard practise in our group to release all research codes as open source~\citep{wang2016three, jiang2016multiphysics, rudraraju2019computational}, the code base of the current work is made available to the wider research community as an open source library~\citep{kunalDendritic2022}.

\subsection{Material properties and non-dimensional quantities in solidification}
{\label{sub:InputParameters}}
We now discuss the specific material properties and non-dimensional parameters used in the modeling of pure metal and binary alloys.  First, we consider the case of dendrite growth in an undercooled melt of pure metal, whose governing equations are given by Equations~\ref{eq:EquationU} and \ref{eq:EquationPHI}. Undercooling of the liquid melt is responsible for driving dendrite growth in this case and there is no external imposition of a temperature gradient or a dendritic growth velocity. This is commonly referred to as free dendritic growth. The set of non-dimensional input parameters used in this model and their numerical values (adopted from Karma and Rappel~\cite{karma1998quantitative}) are summarized in Table~\ref{table:PureMeltParameters}. Here, $\tilde D$ is the thermal diffusivity of heat conduction, $\epsilon_{4}$ is the strength of interface surface energy anisotropy where the subscript indicates the assumed symmetry of the dendritic structure, and in this simulation, we consider a four-fold symmetry. The choice of the length-scale parameter $\lambda_0$, identified as the interface thickness, and the time-scale parameter $\tau_{0}$ set the spatial and temporal resolution of the dendritic growth process. $\xi$ is a coupling parameter representing the ratio of interface thickness and capillary length. 
 

\begin{table}[ht!] 
\centering 
\begin{tabular}{c c c} 
\hline 
Non-dimensional parameter & Numerical value \\  
\hline 
Thermal diffusivity,  $\tilde D$ & 1    \\ 
Anisotropy strength,  $\epsilon_{4}$ &0.05 \\
Characteristic length-scale,  $\lambda_{0}$ & 1  \\
Characteristic time-scale,  $\tau_{0}$ &1 \\
Coupling parameter,  $\xi$ &1.60  \\
\hline 
\end{tabular}
\caption{Input parameters for the phase-field model of solidification of a pure metal. The length-scale parameter $\lambda_0$ and the time-scale parameter $\tau_{0}$ set the spatial and temporal resolution of the dendritic growth process.} \label{table:PureMeltParameters} 
\end{table}

\begin{table}[htp]
\centering 
%
\begin{tabular}{c c c} 
\hline 
Parameter& Numerical value (units) \\ [0.5ex] 
\hline 
Liquid phase diffusivity, $D_{l}$ & $2.4 \times 10^{-9}$ $m^2s^{-1}$ \\ 
Solid phase diffusivity, $D_{s} $ & $1.15 \times 10^{-12}$ $m^2s^{-1}$  \\
Partition coefficient, $k$ & 0.14  \\
Anisotropy strength, $\epsilon_{4}$ &0.05 \\
Gibbs-Thomson coefficient, $\Gamma_g$ & $2.36\times 10^{-7}$ m-K \\
Liquidus slope, $m_{l}$ & -3.5 $K(wt\%)^{-1}$ \\
Initial concentration, $c_{0}$ & 4wt\% \\ 
Thermal gradient, $G$ & 3700 $K m^{-1}$ \\
Pulling velocity, $\bm{\nu}_{p}$ & $1.35 \times 10^{-5}$ $ms^{-1}$ \\
Interface thickness, $\lambda_{0}$ & $1.058 \times 10^{-6}$ m \\
Model constant, $a_{1}$ & 0.8839 \\
Model constant, $a_{2}$ & 0.6267 \\[1ex]
\hline 
\end{tabular}
\caption{Material properties and phase-field model parameters for a binary alloy Al-4wt.\% Cu} 
\label{table:BinaryAlloyProerties} 
\bigskip
\begin{tabular}{c c c} 
\hline 
Input parameters & Expression \\ [0.5ex] 
\hline 
Chemical capillary length,  $d_{0}$& $\frac{k\Gamma}{mC_{0}(k-1)}$  \\[0.5ex]
Ratio of interface thickness to capillary length,  $\Xi$ &$\frac{\lambda_{0}}{d_{0}}$  \\[0.5ex]
Non-dimensional thermal diffusivity,  $\Tilde{D}$& $\frac{D_{l}\tau_{0}}{\lambda_{0}^{2}}$  \\ [0.5ex]
Non-dimensional pulling velocity,  $\Tilde{\bm{\nu}}$& $\frac{\bm{\nu}\tau_{0}}{\lambda_{0}}$  \\[0.5ex]
Characteristic thermal length-scale, $\Tilde{L_{T}}$& $\frac{|m_{l}|(1-k)c_{0}}{kG{\lambda_{0}}}$  \\ [0.5ex] 
Characteristic time-scale,  $\tau_{0}$& $\frac{a_{1}a_{2}\xi \lambda_{0}^{2}}{D_{l}}$  \\[0.5ex]
Coupling parameter,  $\xi$&$a_{1}\Xi$  \\[1ex]
\hline 
\end{tabular}
\caption{Input parameters for the phase-field model of solidification of a binary metallic alloy. The numerical values of the input parameters are computed using the expressions shown here that are in terms of the properties listed in Table~\ref{table:BinaryAlloyProerties}.}
\label{table:PhaseFieldAlloy}
\end{table}

Next we discuss the case of dendritic growth in binary alloys. As mentioned earlier, chemically, alloys are distinct from pure melts due to the presence of solute atoms of one or more alloying materials. During the solidification of alloys, solute diffusion plays a critical role. The coupled governing equations of  solute diffusion and phase-field are given by the Equations~\ref{eqn:PF_Alloy_UC}-\ref{eqn:Alloy_Phi_MODIFIED}. The solid-liquid interface is constrained with the imposition of the temperature gradient,  $G$, and the pulling velocity, $\bm{\nu}_{p}$. We solve for two primal fields - composition undercooling, $u_{c}$, and the phase-field order parameter, $\phi$. The input parameters used in the numerical model are given in Tables~\ref{table:BinaryAlloyProerties}-\ref{table:PhaseFieldAlloy}. The parameters listed in the Table~\ref{table:BinaryAlloyProerties} are the required material properties for modeling dendritic growth in a binary alloy (adopted from Zhu \textit{et al.}\cite{zhu2019phase}): mass diffusivities $D_{l}$ and $D_{s}$ in the solid and liquid phase, Gibbs-Thomson coefficient $\Gamma_g$, thermal gradient, $G$ and pulling velocity $\nu_{p}$. Since here we only model a generic dilute binary alloy, the temperature and composition phase diagram of a physical alloy is approximated by a linear relationship. This approximation is possible if we assume the alloy to be a dilute binary alloy. This linearization of a temperature-composition relationship results in two more parameters: the partition coefficient $k$ and the slope of the liquidus line $m_{l}$. An optimal choice of the phase-field interface thickness is made keeping in mind the computational tractability of the final numerical model. Table~\ref{table:PhaseFieldAlloy} lists the expressions used to calculate the input parameters relevant to the model.   

\subsection{Primary geometric features and shape variations in dendrite morphologies}{\label{sub:ShapeSchematics}}
Dendrites occur in a great variety of morphologies and orientations. One of most common, yet fascinating, occurrences of dendritic growth is in the formation of snowflakes that are dendritic structures with a distinct six-fold geometric symmetry.  
In more practical applications, as in the solidification of metallic alloys, dendritic shapes can be broadly classified as ``columnar'' or ``equiaxed''. Columnar dendrites grow into a positive temperature gradient, and the rate of growth is controlled by the constitutional undercooling of the liquid melt ahead of the dendritic tip.  On the other hand, growth of an equiaxed dendrites in pure melt is controlled by the temperature undercooling of the surrounding liquid melt. 
\par
\begin{figure}[ht!]
  \centering
        \includegraphics[width=0.65\textwidth]{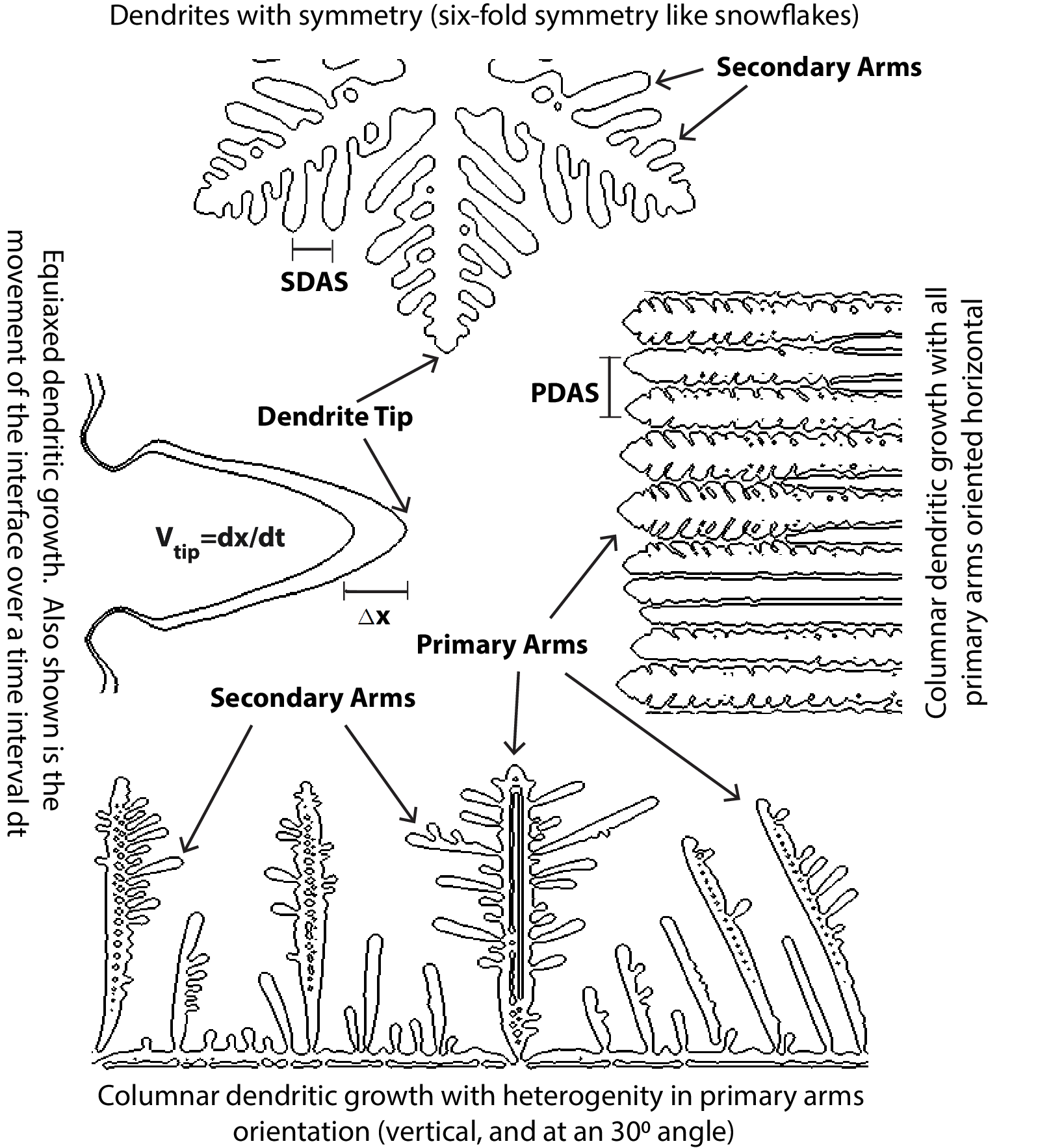}
  \caption{Schematic showing important geometric features used to quantify dendritic morphologies, and some typical variations of shapes and orientations in columnar and equiaxed dendrites.  Shown are the primary arms and primary dendrite arm spacing (PDAS), secondary arms and secondary dendrite arm spacing (SDAS), and the movement of the dendritic tip by a distance $\Delta$x over a time interval $\Delta$t. Dendritic tip velocity is given by $\bm{\nu}_\text{tip}$. All these shapes were generated using the solidification numerical models presented in this work.}
  \label{fig:DendriteShapeSchematic}
\end{figure}

The schematic shown in Figure~\ref{fig:DendriteShapeSchematic} depicts these possible variations in dendritic shapes, and also identifies important geometric features that are traditionally used by the Materials Science and Metallurgy community to quantify dendritic morphologies. At the top end of the schematic, shown is a equiaxed dendrite with a six-fold symmetry (one-half of the symmetric dendrite shown) that is typical of snowflakes and in certain alloy microstructures.  At the bottom, shown are columnar dendrites with their primary arms oriented parallel and inclined at a $30^{\circ}$ angle to the vertical axis. These primary arms develop secondary arms as they evolve, and the separation of the primary arms (referred to as the primary dendrite arm spacing or PDAS) and secondary arms (referred to as the secondary dendrite arm spacing or SDAS) are important geometric characteristics of a given dendrite shape. Further, towards the right side of the schematic is a columnar dendrite colony growing homogeneously oriented along the horizontal axis, and shown towards the left side of the schematic are two time instances (separated by an interval $\Delta$t) of an equiaxed dendrite's interface evolution and the movement of the dendritic tip over this time interval is $\Delta\bm{x}$.  The dendritic tip velocity, a very important measure of dendritic kinetics, is then given by $\bm{{\nu}}_\text{tip}=\frac{d\bm{x}}{dt}$. It is to be noted that all these morphologies were produced by the authors of this manuscript in simulations of various solidification conditions using the numerical framework presented in this work, and are but one demonstration of the capabilities of this framework in modeling complex dendritic morphologies. 

\subsection{Single equiaxed dendrite evolution}{\label{sub:SingleEquiaxed}}
Growth of a single equiaxed dendrite in an undercooled melt of pure metal is a classical solidification problem, and much of traditional analysis of dendritic growth is based on reduced order models for this problem. The phase-field governing equations are given by Equations~\ref{eq:WeakFormU}-\ref{eq:WeakFormPhi}, and the primal fields solved for are the temperature undercooling $u$ and the order parameter $\phi$. These equations are solved using both a $C^0$ basis and a $C^1$ basis as part of the study of convergence with respect to variations in the spatial discretization. Input parameters used in this model have been listed in Section~\ref{sub:InputParameters}. Zero flux Boundary Conditions are considered on all the boundaries for the fields $u$ and $\phi$, and the Initial Condition is a uniformly undercooled liquid melt with $\Delta=-0.75$ everywhere. The phase-field order parameter is set to $\phi=-1.0$ everywhere, except a very small circular region at the center representing a solid seed formed due to nucleation with $\phi=1.0$ in the seed.  Recollect that, in all the models, $\phi=1.0$ represents the bulk solid and $\phi=-1.0$ represents the bulk liquid region and any value in between represents the diffuse interface between these two regions. The non-dimensional size of the domain is  500x500 and it is uniformly discretized into elements with a length measure of  $\Delta x =2.0$. The non-dimensional time step for the implicit time stepping scheme considered in this case is $\Delta t =0.04$. 

\begin{figure}[ht!]
  \centering
  \subfloat[\large Temperature contour at t=840\label{fig:Temperature_SED}]{
        \includegraphics[width=0.49\textwidth]{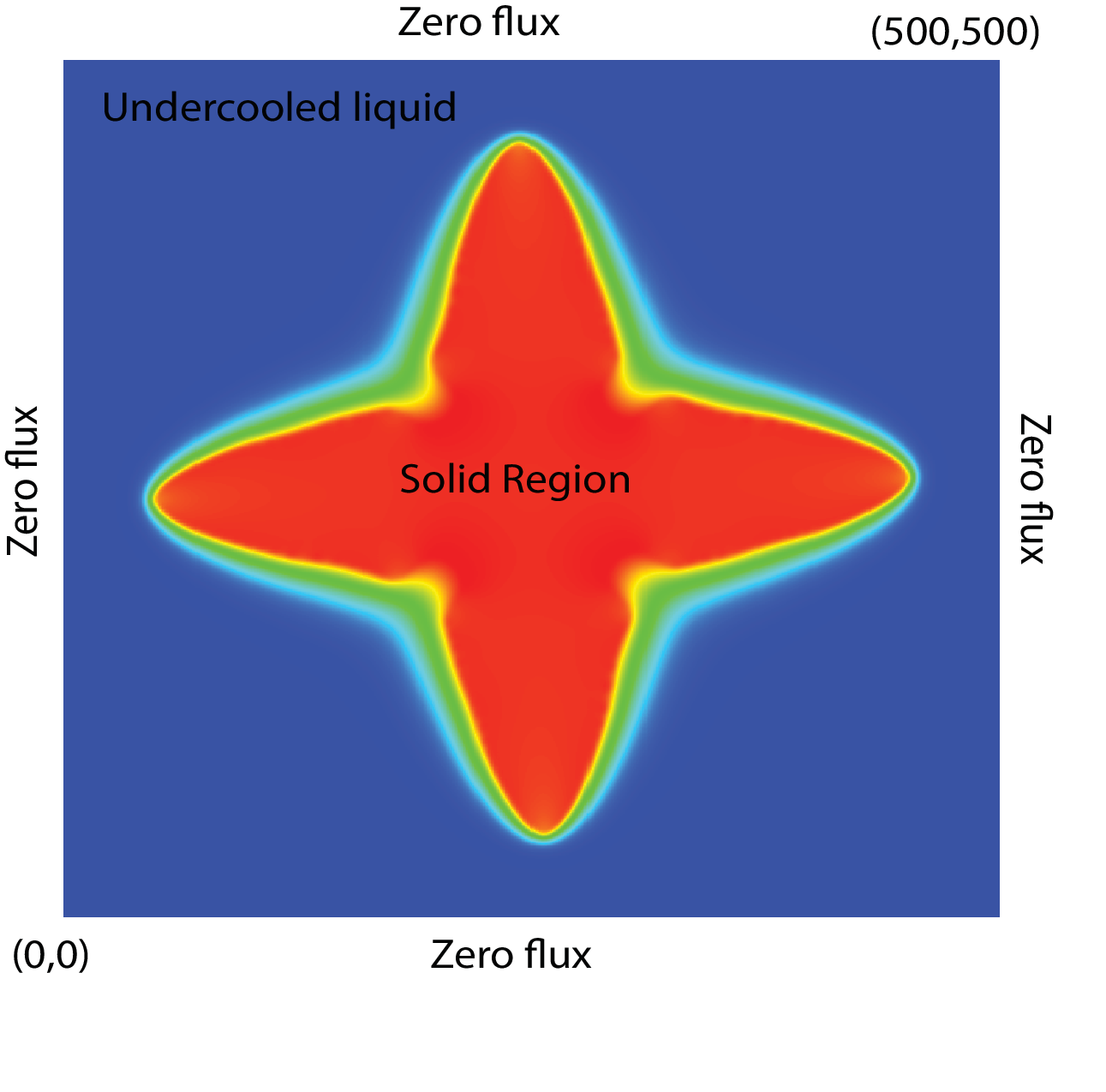}
  }
   \subfloat[\large Time evolution of the interface\label{fig:Phi_SED}]{
        \includegraphics[width=0.49\textwidth]{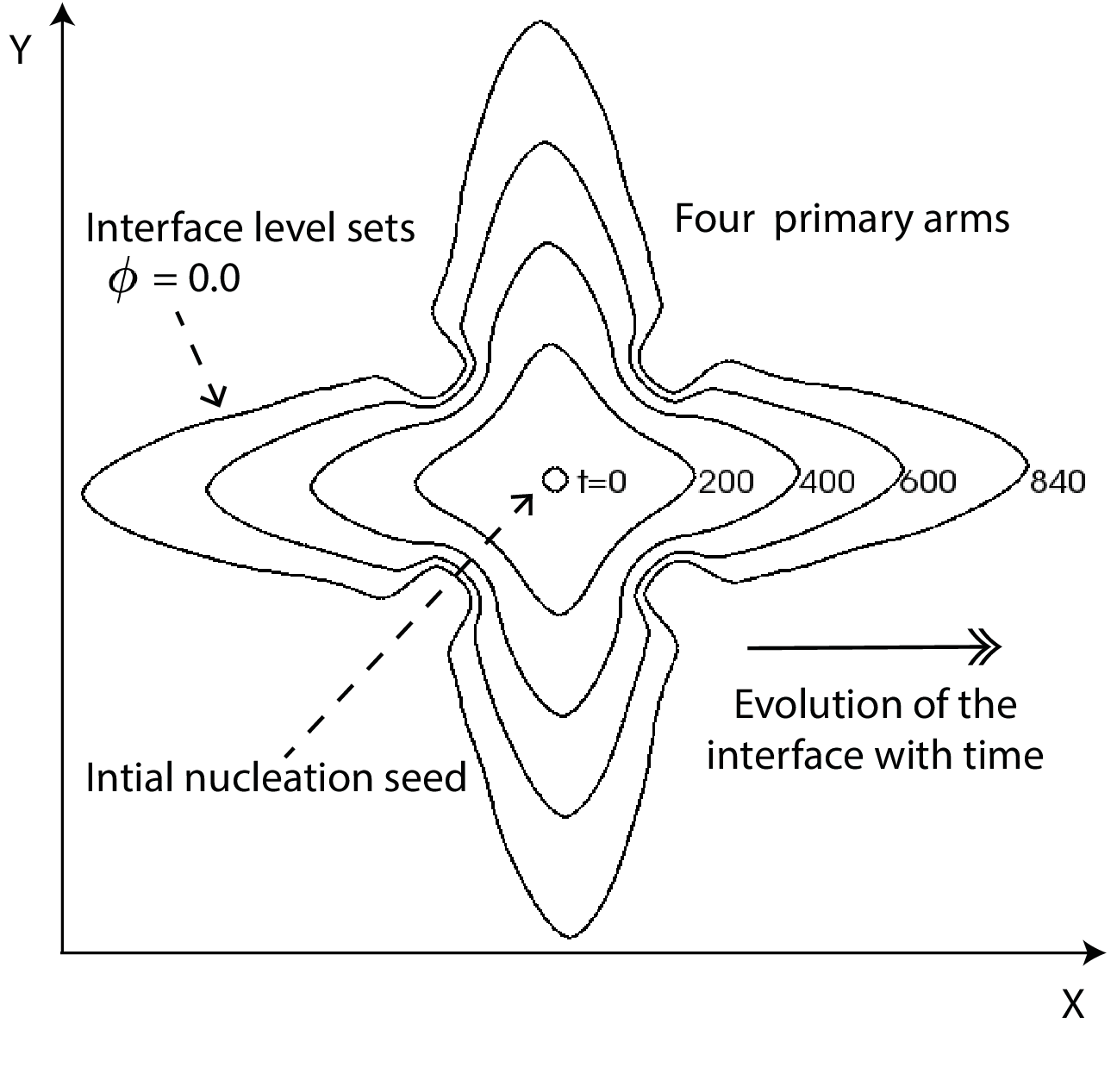}
  }
  \caption{ Growth of a single equiaxed dendrite in an undercooled melt of pure metal. Shown are the (\ref{fig:Temperature_SED}) Temperature contour $(u)$ at t = 840, and (\ref{fig:Phi_SED}) Temporal evolution of the equiaxed dendrite interface, delineated by the level set $\phi=0$, at different time instances.}
  \label{fig:SED}
\end{figure}

The growth of dendrites out of the nucleation point at the center is dictated by the degree of undercooling of the liquid melt.  As solidification begins,  the surrounding undercooled liquid melt is at lower temperature than the solid seed.  As a result, the solid seed diffuses heat into the liquid region. Theoretically, the solid-liquid interface of an undercooled pure melt is always unstable, so small local numerical perturbations on circular solid seed interface begin to grow.  Due to the assumption of strong anisotropy (four-fold symmetry) of the interface, perturbations to the interface grow rapidly along the favorable orientations as compared to the non-favorable orientations. Figure~\ref{fig:Temperature_SED} represents $u$ at a particular time, t = 840. Diffusion takes place over the entire stretch of the solid and liquid regions.  From t = 0 to 840, $u$ changes from the initial value of -0.75 as per the solution of the heat diffusion equation. It should be noted that the bulk of the liquid region and the solid dendrite are at a uniform temperature, and the temperature variations are mostly localized to the vicinity of the interface. The bulk liquid melt region shown in blue is at a low temperature and the bulk solid dendrite in red is at a high temperature, as seen in Figure~\ref{fig:Temperature_SED}. The level set of $\phi=0$ is identified as the distinct solid-liquid interface and the evolution of this interface is plotted in Figure~\ref{fig:Phi_SED} as a series of $\phi=0$ level sets at different time instances. As can be seen, at t=0, the circular seed at the center of the domain, over time, evolves into a four-fold symmetric interface. At time, t = 840, we see a fully developed equiaxed dendrite shape. We also simulated this Initial Boundary Value Problem on $C^1$-continuous mesh and found identical spatiotemporal growth of the equiaxed dendrite. 

\subsubsection{Numerical convergence studies}{\label{sub:OptimalConvergence}}
In this section, we look at the accuracy of the interface propagation numerical results for the single equiaxed dendrite case presented above. We do this by comparing the velocity of propagation of the dendritic tip with known analytical solutions of tip velocity for this problem, and also study the convergence of the tip velocity with mesh size and increasing continuity of the numerical basis.     

\paragraph{Velocity of the dendritic tip:}{\label{sub:TipVelocity}}
As can be seen in Figure~\ref{fig:SED}, the phase-field model of solidification of a pure metal successfully produces the freely growing dendrite shape and four-fold symmetry as expected. However, predicting quantifiable features of the dendritic evolution would provide a better validation of the proposed numerical framework. In this context, the dendritic tip propagation velocity (see Figure~\ref{fig:DendriteShapeSchematic}) for single equiaxed dendrites has been widely used in analytical studies of dendritic propagation. In the literature, implicit expressions for radius and velocity of the dendritic tip have been obtained using semi-analytical approaches, and they provide estimates of $\bm{\nu}_{tip}$. One of the more popular semi-analytical approaches to obtain the tip velocity employs a Green's function (GF) approach. The Green's function is an integral kernel that is widely used to solve of certain classes of differential equations. We will not discuss this approach here, but only use the dendritic tip velocity values predicted by this approach. These values are denoted by $\bm{\nu}_{tip}(GF)$. Interested readers are referred to Karma and Rappel\cite{karma1998quantitative} for details of the GF approach.  
\par
The model parameters used in this study are $\Delta=-0.65$, $\epsilon_{4}=0.05$, $\tilde D=1$, $d_{0}/\lambda_{0}=0.544$, $\lambda_{0}=1$, $\tau_{0}=1$, $\xi=1.596$. The numerical domain size and time step size, in dimensionless units, are 500x500 and $\Delta t=0.04$, respectively.  \textcolor{black}{We chose an implicit (backward Euler scheme) time-stepping method for this study, and the time step size, $\Delta t$, was chosen considering stable quadratic convergence of the relative residual error for the implicit time stepping scheme. Further, this time step size is comparable to the values reported in the finite difference method literature for this problem \cite{karma1998quantitative}}. For these values, the Green's function approach predicts an equilibrium tip velocity,  $\bm{\nu}_{tip}(GF)=0.0469$. Freely growing equiaxed dendrite in an infinite domain, under solidification conditions held constant over time, attains a fixed equilibrium tip velocity. In this study, we consider the effect of two discretization parameters, namely: (1) Comparison of the numerical dendritic tip velocity ($\bm{\nu}_{tip}$) to the equilibrium tip velocity obtained from Green's function approach $\bm{\nu}_{tip}(GF)$ (Figure~\ref{fig:VtipCvsTime1}), and (2)  Convergence of the tip velocity ($\bm{\nu}_{tip}$) with mesh refinement, and its dependence on the order of continuity of the basis (Figure~\ref{fig:VtipCvsDeltaX}). 
\par 
The method of estimating the dendrite tip velocity from the numerical simulations using the computational frameworks described above is described briefly. Simulations were performed using both a $C^0$-continuous basis and a $C^1$-continuous basis, to study the effect of basis continuity when resolving fine dendritic morphologies. The dendritic tip velocity measurement is done manually as a post-processing of the phase-field order parameter contours of the numerical simulations. \textcolor{black}{At each time instance, the zero level set ($\phi=0$) is plotted using a visualization tool for FEM output, Visit~\cite{Visit}. The tip of the zero level set of the order parameter is manually tracked, and the velocity of tip propagation at each time instance is computed using a simple finite difference estimate of the derivative $\bm{\nu}_{tip}=\frac{d\bm{x}}{dt}$.  The criteria for choosing the zero level set to represent the interface is a standard practice in the dendritic modeling literature \cite{boettinger2002phase, gibou2003level}. As can be seen from Figure~\ref{fig:VtipCvsTime1}, after a non-dimensional time of 1500, for all the cases we considered, the estimated tip velocity at any time instant is within 10\% of the theoretical tip velocity, with 10\% being the upper bound for the coarsest $C^0$-continuous mesh. Thus, for the error estimates, we consider the tip velocity at t = 1500. At time t=2000 the dendrite approaches the edge of the computational domain}. Two sources of error can potentially cause minor variations to the computed velocity values. The first source of error arises from not picking the exact tip location in subsequent time snapshots as the dendrite interface shape is evolving with time. The second source of error arises due to the use of only $C^0$ interpolation of the nodal solutions inside the elements by the visualization tool. This error is significant for $C^1$-basis, as the results are $C^1$-continuous but the visualization of the results inside the elements is only $C^0$-continuous.  As can be seen in Figure~\ref{fig:VtipCvsTime1}, $\bm{\nu}_{tip}$ starts to approach the Green's function value, $\bm{\nu}_{tip}(GF)$, across all the simulations, but shows better convergence to the Green's function value for the $C^1$-continuous basis. Further, we observe the percentage error of $\bm{\nu}_{tip}$, given by Error$=\frac{|\bm{\nu}_{tip}-\bm{\nu}_{tip}(GF)|}{\bm{\nu}_{tip}(GF)}$, consistently reduces with a decrease in the mesh size across both $C^0$-continuous and $C^1$-continuous basis, as shown in Figures~\ref{fig:VtipCvsDeltaX}(a) and \ref{fig:VtipCvsDeltaX}(b). This is an  important numerical convergence result, as to the best of our knowledge, this is the first time convergence rates have been demonstrated for non-trivial dendritic shape evolution with respect to mesh refinement. The convergence rates are demonstrated in the plot in Figures~\ref{fig:VtipCvsDeltaX}(a), where the rate of convergence (slope of the percentage Error vs mesh size plot) is close to two for a $C^0$-continuous basis, and close to three for a $C^1$-continuous basis. While obtaining theoretical optimal converge rates with respect to h-refinement of the discretization for this coupled system of parabolic partial differential equations is very challenging, one can observe that the convergence rates shown are comparable to the optimal convergence rates of classical heat conduction and mass diffusion partial differential equations~\citep{bieterman1982finite}. This completes the estimation of error in $\bm{\nu}_{tip}$, which is a measure of localized error in the order parameter field, as the dendritic interface is given by the zero level set of the order parameter field. Now we look at estimation of errors over the entire domain for both the temperature field and the order parameter field.  
 
\pgfplotstableread{Textfiles/VtipvsTime.txt}{\VtipCvsTime} 
\pgfplotstableread{Textfiles/VtipvsTimecone.txt}{\VtipCvsTimecone}
\pgfplotstableread{Textfiles/VtipvsTimectwo.txt}{\VtipCvsTimectwo}

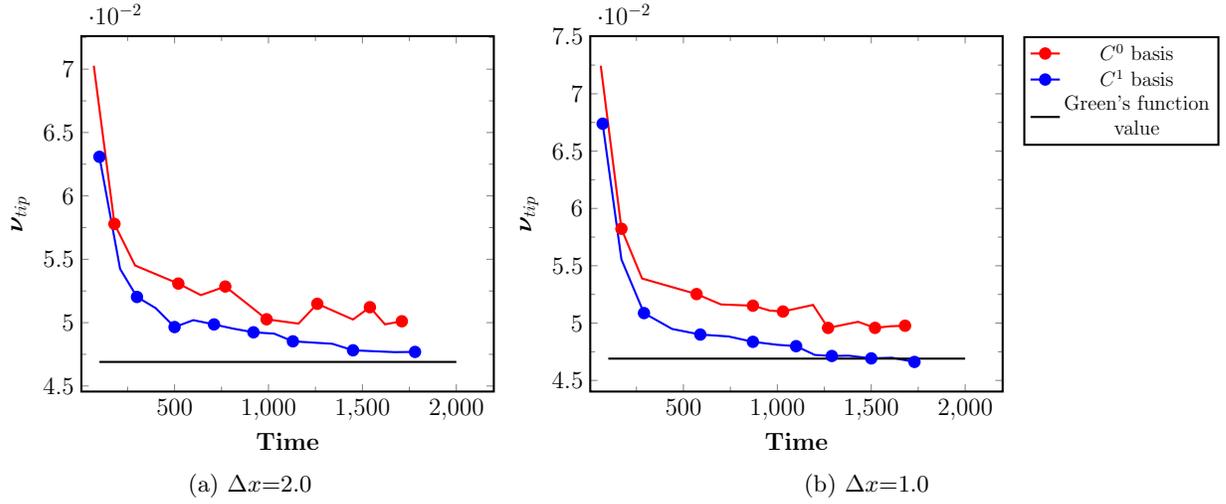
\begin{figure}[ht!]
  \centering
  \subfloat[$\Delta x$=2.0\label{fig:VtipVsTimeX2.0_C0C1TH}]{
    \begin{tikzpicture}[scale=0.65]
      \begin{axis}[ticklabel style = {font=\Large },very thick,minor tick num=1,xlabel={\Large \textbf{Time}},ylabel={\Large $\boldsymbol{\nu}_{tip}$}, x unit=, y unit=, mark repeat={2}, scale only axis=true, xmin=4.0] 
        \addplot [red, mark=*, mark size=3] table [x={t_one}, y={v_one}] {\VtipCvsTime};
         \addplot [blue, mark=*, mark size=3] table [x={t_one}, y={v_one}] {\VtipCvsTimecone};
        \addplot [black, no markers] coordinates {(100,0.0469) (2000,0.0469)};
      \end{axis}
    \end{tikzpicture}
  }  
\subfloat[$\Delta x$=1.0\label{fig:VtipVsTimeX1.6_C0C1TH}]{
    \begin{tikzpicture}[scale=0.65]
      \begin{axis}[ticklabel style = {font=\Large },very thick, minor tick num=1,xlabel={\Large \textbf{Time} },ylabel={\Large  $\boldsymbol{\nu}_{tip}$}, x unit=, y unit=, mark repeat={2}, scale only axis=true,xmin=4.0, legend style={cells={align=center}}] 
        \addplot [red, mark=*, mark size=3] table [x={t_two}, y={v_two}] {\VtipCvsTime};
           \addplot [blue, mark=*, mark size=3] table [x={t_two}, y={v_two}] {\VtipCvsTimecone};
        \addplot [black, no markers] coordinates {(100,0.0469) (2000,0.0469)};
       \legend{{\large $C^0$ basis}, {\large $C^1$ basis},{\large Green's function \\ \large value}}
      \end{axis}
    \end{tikzpicture}
  }
  \caption{ Comparison of the numerical and analytical (Green's function approach) dendritic tip velocities over time during growth of a single equiaxed dendrite for two different mesh sizes ($\Delta x$), and for both $C^0$-continuous and $C^1$-continuous basis.}
  \label{fig:VtipCvsTime1}
\end{figure}

\pgfplotstableread{Textfiles/TipError.txt}{\TipError}
\begin{figure}[ht!]
  \centering
   \begin{minipage}[c]{0.45 \linewidth}
    \centering
    \begin{tikzpicture}[scale=0.85]
      \begin{loglogaxis}[ticklabel style = {font=\large },very thick, minor tick num=1,xlabel={\large ${\Delta x}$},ylabel={\large Error}, x unit=, y unit=, mark repeat={0}] 
        \addplot [red,only marks, mark=o, mark size=3] table [x={C0h}, y={C0Error}] {\TipError};
        \addplot [red, mark=none,line width=1pt] table [x={C0h}, y={FitC0}] {\TipError};
         \addplot [blue,only marks, mark=o,  mark size=3] table [x={C1h}, y={C1Error}] {\TipError};
         \addplot [blue, mark=none,line width=1pt] table [x={C1h}, y={FitC1}] {\TipError};
        \addplot [black, no markers] coordinates {(1.2,2.0) (1.2,3.6)};
       \addplot [black, no markers] coordinates {(0.92,2.0) (1.2,2.0)};
        \addplot [black, no markers] coordinates {(1.8,1.2) (1.8,2.4)};
       \addplot [black, no markers] coordinates {(1.44,1.2) (1.8,1.2)};
      
        \node[below] at (axis cs:1.23,3.1){\large $\textbf{2}$};
        \node[below] at (axis cs:1.06,2.0){\large $\textbf{1}$};
        \node[below] at (axis cs:1.85,1.82){\large $\textbf{3}$};
        \node[below] at (axis cs:1.62,1.2){\large $\textbf{1}$};
        
        
        \node[fill =white] at (axis cs:1.9,0.42) {$y_{0}=2.04x_{0}+0.35$};
        \node[fill =white] at (axis cs:1.9,0.22) {$y_{1}=3.72x_{1}-0.48$};
        \node[fill =white] at (axis cs:1.9,0.12) {$y_{2}=4.11x_{2}-0.51$};
         \legend{{}, {\large $C^{0}$ basis}, {},{\large $C^1$ basis},{},{}}
      \end{loglogaxis}
    \end{tikzpicture}
   \caption*{(a) Convergence of error in dendritic tip velocity with mesh refinement.}
  \label{fig:ErrorinTipVelocity}
  \end{minipage}
  \hfill
   \begin{minipage}[c]{0.45 \linewidth}
   \centering
  \centering
   \begin{tikzpicture}[scale=0.88]
    \begin{axis}[ticklabel style = {font=\large },very thick,
    ybar,bar width=30,
    enlargelimits=0.45,
    ylabel= \large Error,
    symbolic x coords={$C^0$ basis, $C^1$ basis},
    xtick=data,
    legend style={at={(0.5,0.97)}, 
			anchor=north,legend columns=-1},
		domain=2:8
    nodes near coords, 
	nodes near coords align={vertical}]
    \addplot coordinates {($C^0$ basis,8.11) ($C^1$ basis,1.58)};
		\legend{\large $\Delta x$ = 1.0} 
  \end{axis}
    \end{tikzpicture}
  \caption*{(b) Comparison of error in dendritic tip velocity for a mesh size of $\Delta x$ = 1.0}
  \label{fig:BarPlot}
   \end{minipage}
    \caption{Convergence of the percentage error in dendritic tip velocity, given by Error=$\frac{|\bm{\nu}_{tip}-\bm{\nu}_{tip}(GF)|}{\bm{\nu}_{tip}(GF)}$, with mesh refinement and its dependence on the continuity of the basis. \textcolor{black}{The $C^0$ basis is chosen to be first order (Linear), and the $C^1$ basis is chosen to be second order (Quadratic). Subplot (a) shows rate of error convergence ($\epsilon = \mathcal{O}(\Delta x)^{order+1}$) with respect to mesh refinement in a log-log plot}. Subplot (b) shows the difference in percentage errors between a $C^0$-continuous basis and a $C^1$-continuous basis for one mesh size.}
    \label{fig:VtipCvsDeltaX}
\end{figure}

\paragraph{\textcolor{black}{Error convergence of the temperature and order parameter fields:}}{\label{sub:SpatialDiscretization}}
For this study, we consider three different mesh refinements of a $C^0$-continuous basis. For each of these refinements, the computational domain size is 160x160, chosen time step is $\Delta t = 0.12$, and the other input parameters are as given in Table ~\ref{table:PureMeltParameters}.  \textcolor{black}{We chose an implicit (backward Euler scheme) time-stepping method for this study, and the time step size, $\Delta t$, was chosen considering stable quadratic convergence of the relative residual error for the implicit time stepping scheme. Further, we also verified that a time-step size smaller than the value considered does not yield any significant improvement in the error computations}. As the Initial Condition, a small nucleation seed was considered, and this seed grows into the surrounding undercooled melt. For the problem of interest, there are no known analytical solutions for the fields $u$ or $\phi$. As mentioned in the earlier subsection, the only analytical result available is for the equilibrium dendritic tip velocity. Thus, we follow a computational approach to obtain the reference solution of both the primal fields. A separate problem is solved on a very fine mesh ($\Delta x =$ 0.2) using the popular open source phase field library PRISMS-PF~\citep{dewitt2020prisms}, and the primal field solutions for this very fine mesh are considered the reference solutions. It is to be noted that one of the authors (S.R.) in this work was the lead developer of the PRISMS-PF library, and this library is being successfully used by the wider Materials Science community to solve numerous phase-field problems.  \textcolor{black}{The mesh size of $\Delta x =$ 0.2 for the reference solution was arrived at by performing a convergence study of the total free energy with respect to the mesh refinement, and noting the mesh size below which the rate of change of free energy with respect to mesh refinement nearly plateaued.} 
\par
Using this numerical reference solution, error in the temperature and order parameter fields of the solutions obtained for the numerical formulations presented in this work, is calculated by comparing the fields and their gradients using two types of norms relevant to the problem, and for three different mesh refinements at one instance of time (time = 450) when the dendritic shape is well established. The two error norms used in this analysis are the $L_{2}$-norm and the $H_{1}$-norm, and these are defined as follows: 
\begin{subequations}
\begin{equation}
|u-\hat u|_{L_2} =\sqrt{ \int_{\Omega}|u_{h}-\hat u|^2~dV} 
\end{equation}
\begin{equation}
|u-\hat u|_{H_{1}}=\sqrt{\int_{\Omega} \Big(|u_{h}-\hat u|^2 +|\nabla u_{h}-\nabla \hat u|^2\Big)~dV} 
\end{equation}
\end{subequations}
\noindent Here, $u_{h}$ is the finite element solution under consideration and $\hat u$ is the reference numerical solution. The integral is defined over the entire numerical domain, so the estimated error, unlike the previous error in tip velocity, is not localized to any particular region in the domain. The convergence of this error in these two norms with mesh refinement is shown in Figure~\ref{fig:LogErrosVsLogh} on a log-log scale. As can be seen from these plots, we observe good convergence in both the norms for both the primal fields, with the rates of convergence being in the vicinity of 1.5 to 2.0.   \textcolor{black}{In this work, we opted to study the error convergence for each of the primal fields separately. Given that this is a coupled problem in two primal fields, it is possible to construct a single $L_{2}$-norm and $H_{1}$-norm that accounts for both the fields. For example, $\sqrt{ \int_{\Omega} ( |u_{h}-\hat u|^2 + |\phi_{h}-\hat \phi|^2)~dV}$ for a unified $L_2$-norm of the error in both the fields.} We are not aware of any theoretical optimal error estimates in any of the relevant error norms for this problem, so this numerical strategy to obtain error estimates is seen as filling this vital gap by providing numerical convergence studies, thus permitting the numerical error convergence analyses of these governing equations.

\pgfplotstableread{Textfiles/Temperature1.txt}{\Temperature}
\pgfplotstableread{Textfiles/PhaseField1.txt}{\PhaseField}
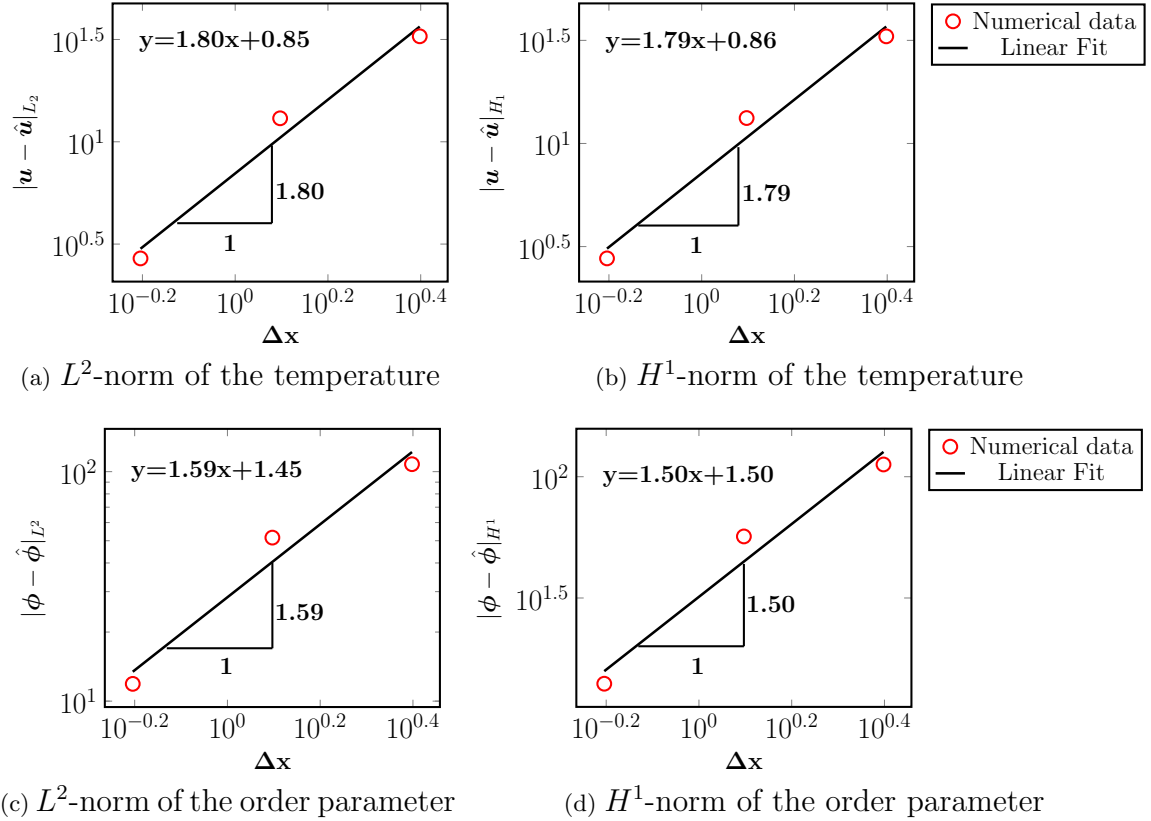
\begin{figure}[ht!]
  \centering
  \subfloat[\large $L^2$-norm of the temperature\label{fig:UTL2}]{
    \begin{tikzpicture}[scale=0.65]
      \begin{loglogaxis}[ticklabel style = {font=\LARGE }, minor tick num=1,xlabel={\Large $\mathbf{\Delta x}$ },ylabel={\Large  $|\boldsymbol{u}-\hat{\boldsymbol{u}}|_{L_{2}}$}, x unit=, y unit=, mark repeat={0},very thick] 
        \addplot [red,only marks, mark=o, mark size=4] table [x={dx}, y={L2}] {\Temperature};
        \addplot [black, mark=none,line width=1.5pt] table [x={dx}, y={FitL2}] {\Temperature};
        \addplot [black, no markers] coordinates {(1.2,4.0) (1.2,9.6)};
       \addplot [black, no markers] coordinates {(0.75,4.0) (1.2,4.0)};
        \node[below] at (axis cs:0.975,3.8){\Large $\textbf{1}$};
        \node[below] at (axis cs:1.38,6.8){\Large $\textbf{1.80}$};
        \node[fill =white] at (axis cs:0.95,31) {\Large $\textbf{y=1.80x+0.85}$};
      \end{loglogaxis}
    \end{tikzpicture}
  }
 \subfloat[\large $H^1$-norm of the temperature\label{fig:UTH1}]{
    \begin{tikzpicture}[scale=0.65]
      \begin{loglogaxis}[ticklabel style = {font=\LARGE },very thick,minor tick num=1,xlabel={\Large $\mathbf{\Delta x}$ }, ylabel={\Large  $|\boldsymbol{u}-\hat { \boldsymbol{u}}|_{H_{1}}$},  x unit=, y unit=, mark repeat={0}] 
        \addplot [red,only marks, mark=o,  mark size=4] table [x={dx}, y={H1}] {\Temperature};
        \addplot [black, mark=none,line width=1.5pt] table [x={dx}, y={FitH1}] {\Temperature};
        \addplot [black, no markers] coordinates {(1.2,4.0) (1.2,9.6)};
       \addplot [black, no markers] coordinates {(0.73,4.0) (1.2,4.0)};
        \node[below] at (axis cs:0.975,3.8){\Large $\textbf{1}$};
        \node[below] at (axis cs:1.38,6.8){\Large $\textbf{1.79}$};
        
        \legend{{\Large Numerical data}, {\Large Linear Fit}}
        \node[fill =white] at (axis cs:0.95,31) {\Large $\textbf{y=1.79x+0.86}$};
      \end{loglogaxis}
    \end{tikzpicture}
  }
  
  \subfloat[\large $L^2$-norm of the order parameter\label{fig:UPHiL2}]{
    \begin{tikzpicture}[scale=0.65]
      \begin{loglogaxis}[ticklabel style = {font=\LARGE },very thick,minor tick num=1,xlabel={\Large $\mathbf{\Delta x}$ },ylabel={\Large $|\bm{\phi}-\hat{\bm{\phi}}|_{L^2}$}, x unit=, y unit=, mark repeat={0}] 
        \addplot [red,only marks, mark=o,  mark size=4] table [x={dx}, y={L2}] {\PhaseField};
        \addplot [black, mark=none,line width=1.5pt] table [x={dx}, y={FitL2}] {\PhaseField};
        \addplot [black, no markers] coordinates {(1.25,17.0) (1.25,40.0)};
       \addplot [black, no markers] coordinates {(0.74,17.0) (1.25,17.0)};
        \node[below] at (axis cs:0.995,16.5){\Large $\textbf{1}$};
        \node[below] at (axis cs:1.43,28.5){\Large $\textbf{1.59}$};
        
         \node[fill =white] at (axis cs:0.95,100) {\Large $\textbf{y=1.59x+1.45}$};
      \end{loglogaxis}
    \end{tikzpicture}
  }
 \subfloat[\large $H^1$-norm of the order parameter\label{fig:UPHIH1}]{
    \begin{tikzpicture}[scale=0.65]
      \begin{loglogaxis}[ticklabel style = {font=\LARGE },very thick,minor tick num=1,xlabel={ \Large $\mathbf{\Delta x}$  },ylabel={\Large $|\bm{\phi}-\hat{\bm{\phi}}|_{H^1}$}, x unit=, y unit=, mark repeat={0}] 
        \addplot [red,only marks, mark=o,  mark size=4] table [x={dx}, y={H1}] {\PhaseField};
        \addplot [black, mark=none,line width=1.5pt] table [x={dx}, y={FitH1}] {\PhaseField};
         \addplot [black, no markers] coordinates {(1.25,20.0) (1.25,43.7)};
       \addplot [black, no markers] coordinates {(0.74,20.0) (1.25,20.0)};
        \node[below] at (axis cs:0.995,19.1){\Large $\textbf{1}$};
        \node[below] at (axis cs:1.43,34.5){\Large $\textbf{1.50}$};
        
        \legend{{\Large  Numerical data}, {\Large Linear Fit}}
        \node[fill =white] at (axis cs:0.95,100) {\Large $\textbf{y=1.50x+1.50}$};
      \end{loglogaxis}
    \end{tikzpicture}
    }
\caption{Error estimates and order of convergence for the primal fields $(u,\phi)$ shown on log-log plots. The error norms used in this study are the $L_{2}$-norm and the $H_{1}$-norm of $(u,\phi)$. \textcolor{black}{A first order $C^0$-continuous basis is used for this study}. The equation of the linear fit shown is noted in each plot.}
  \label{fig:LogErrosVsLogh}
\end{figure}

\subsection{Growth of multiple equiaxed dendrites} 
\begin{figure}[ht!]
  \centering
  \subfloat[\large Composition\label{fig:Composition_ME}]{
       \includegraphics[width=0.44\textwidth]{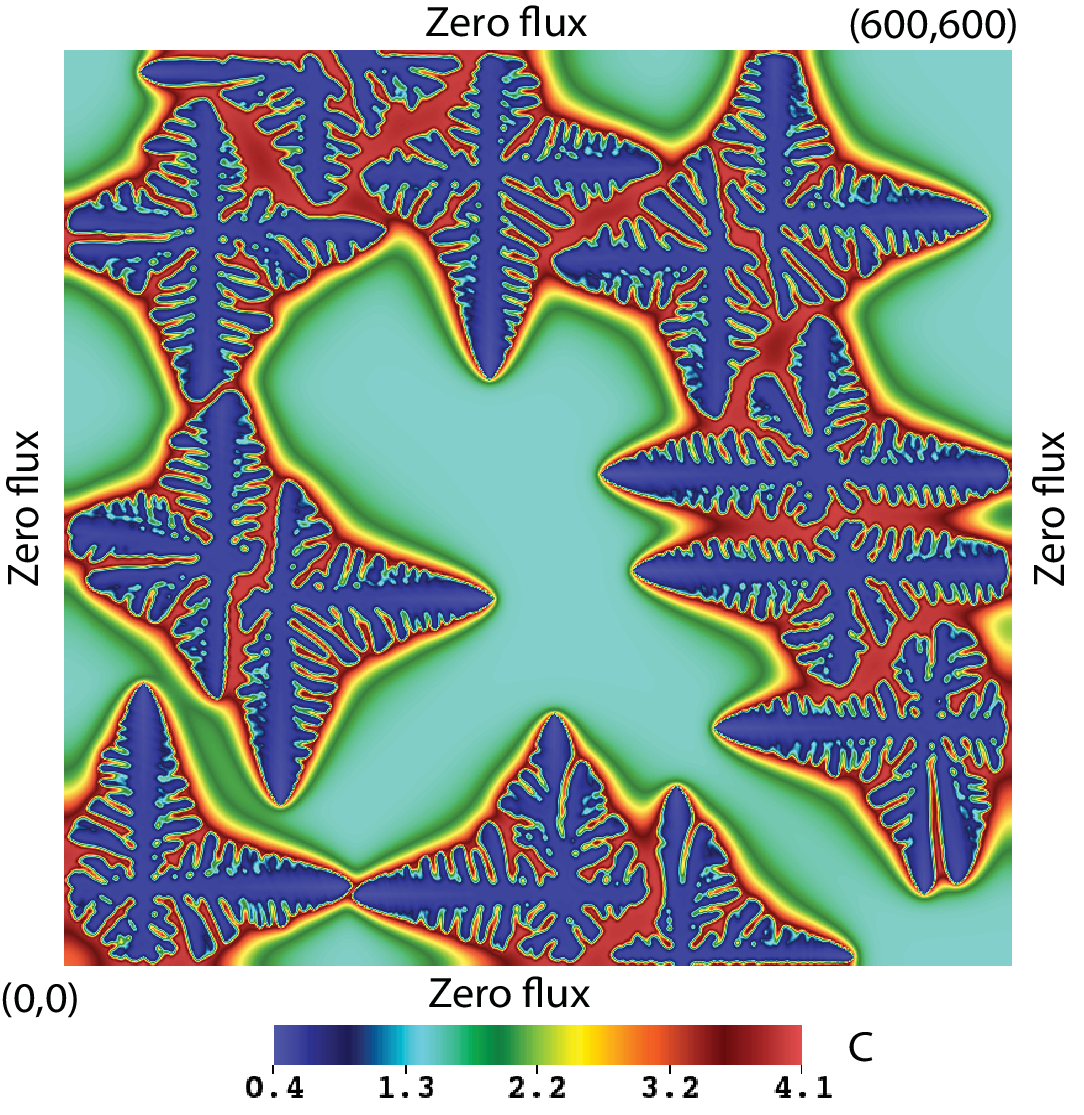}
  }
   \subfloat[\large Zero level sets of the order parameter\label{fig:Phi_ME}]{
        \includegraphics[width=0.54\textwidth]{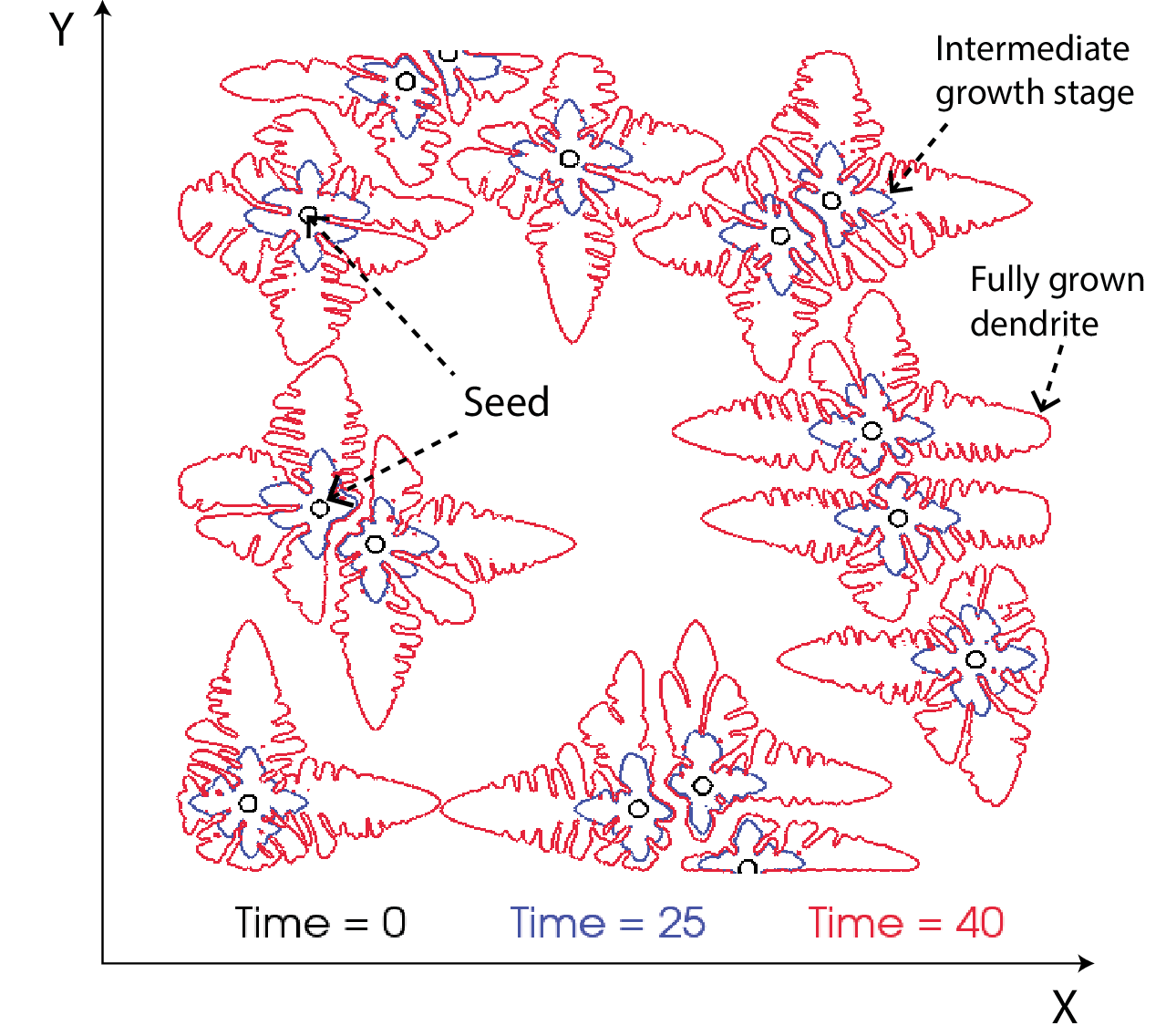}
  }
  \caption{ Growth of multiple equiaxed dendrites. Shown are the (\ref{fig:Composition_ME}) contours of solute composition, c,  at time t = 40, and (\ref{fig:Phi_ME}) Time evolution of multiple equiaxed dendrite interfaces at time t = 0, 25 and 40. The dendrite interfaces are the zero level sets of the order parameter, $\phi=0.0$.}
  \label{fig:EquiaxC1}
\end{figure}

{\label{sub:MultipleEquiaxedDendrite}}
In Section \ref{sub:SingleEquiaxed}, we discussed one particular type of dendritic growth, namely the evolution of a single equiaxed dendrite during the solidification of pure metal. Now we look at the evolution of multiple equiaxed dendrites during the solidification of a more complex binary alloy solution. The input parameters for the model are listed in Table ~\ref{table:BinaryAlloyProerties}. Growth of multi-equiaxed dendrites during solidification of a binary alloy is governed by the Equations~\ref{eq:WeakFormU_Alloy}-\ref{eq:WeakFormPhi_Alloy}. The primal fields in this problem are the non-dimensional composition, $u_{c}$, and the phase-field order parameter, $\phi$. As an Initial Condition for this problem, we consider multiple small nucleation seeds randomly distributed in the numerical domain. The seed locations are completely random, and generated using a random number generator. To ensure free equiaxed dendritic growth, thermal gradient $G$ and pulling velocity $\bm{\nu}_{p}$ is taken to be zero. The numerical domain of size 600x600 is subdivided into uniform elements with a length measure of $\Delta x = 1.2$. The time step size is $\Delta t= 0.005$. The average alloy composition in the entire domain is taken to be $c_{0}=4\%$, which corresponds to $u_{c}=0$. Nucleation seeds of the solid phase are randomly placed such that $\phi=1$ inside of the seed and $\phi=-1$ outside of the seed. Zero flux Boundary Conditions are considered for both the fields. During the metal casting process, many solid nuclei are randomly formed near cavities or defects on the surface of a mold.  These nuclei then start to grow in an equiaxial manner with neighboring nuclei competing with each other. This competition can shunt the growth of dendrites approaching each other. Solution contours of the fields $u_{c}$ and $\phi$ are shown in the Figure~\ref{fig:EquiaxC1}. Here, Subfigure~\ref{fig:Composition_ME} shows the solute composition in the numerical domain at t = 40. Initially, the entire domain is at an alloy concentration of 4\%. Due to negligible mass diffusivity in the solid phase, solute diffusion primarily takes place in the liquid region.  A solute composition of $c_{0}=4\%$ is highest near the interface, and as one moves away from the interface, the solute concentration decreases. This is apparent in all the dendritic interfaces shown in the Subfigure~\ref{fig:Composition_ME}. Solute composition inside the solid is given by $c_{0}=0.56$ and it stays constant due to low negligible diffusivity in the solid region. Subfigure~\ref{fig:Phi_ME} demonstrates the competition between the dendrites growing out of the various seeds, the interaction of the dendritic tips, and the growth of secondary dendrite arms. Material properties of an alloy are often correlated with fine-scale microstructure length scales (at the grain scale and sub-grain scale), and the secondary dendrite arm spacing (SDAS) is one of the important microstructural properties of alloys. 

\subsection{Growth of single columnar dendrites}
\begin{figure}[ht!] 
  \centering
  \subfloat[\large Composition \label{fig:Composition_SSC}]{
        \includegraphics[width=0.45\textwidth]{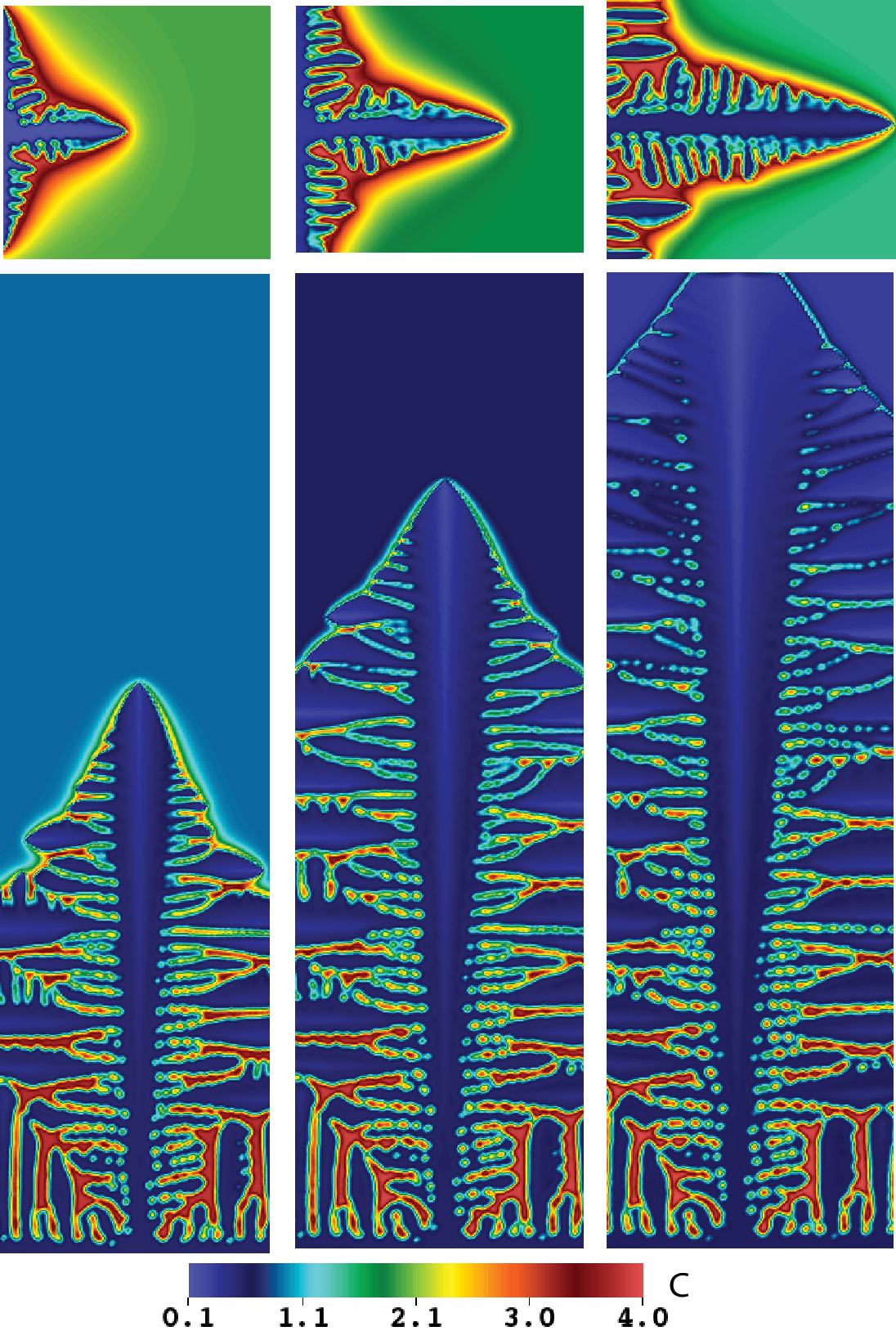}
     }
  \subfloat[\large Zero level sets of the order parameter \label{fig:Phi_SSC}]{
        \includegraphics[width=0.5\textwidth]{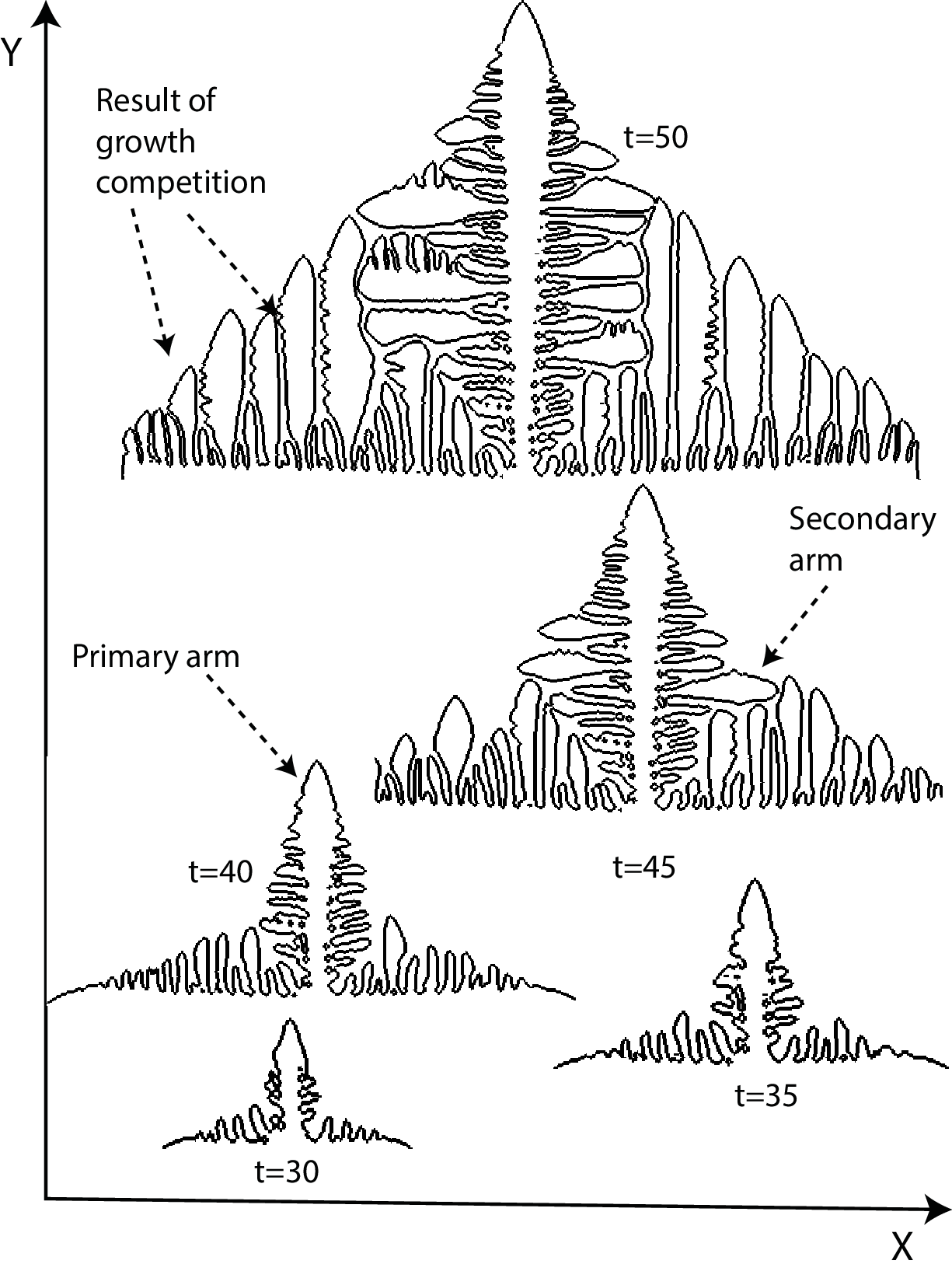}
     }
  \caption{Growth of a single columnar dendrite. Shown are the (\ref{fig:Composition_SSC}) contours of solute composition, c, at different time instances, and (\ref{fig:Phi_SSC}) Time evolution of single columnar dendrite interfaces at times t = 30 to 50.}
  \label{fig:SSColumnar}
\end{figure}

{\label{sub:SingleSeedColumnar}}
In this section, we demonstrate the growth of columnar dendrites, growing from a single seed during the solidification of a binary alloy. The numerical model and input parameters used in this case are identical to the previous Section~\ref{sub:MultipleEquiaxedDendrite}. The notable difference from the model in Section~\ref{sub:MultipleEquiaxedDendrite} to the current model is the imposition of constraints like non-zero thermal gradient $G$ and pulling velocity $\bm{\nu}_{p}$. Initially, we place a single solid seed at the bottom of the numerical domain. The entire domain is at a solute composition of $c_{0}=4\%$ which corresponds to the composition undercooling of $u_{c}=0$. Zero flux Boundary Conditions are enforced on both the primal fields, $u_{c}$ and $\phi$. Figure~\ref{fig:SSColumnar} depicts the evolution of the composition field, and the dendritic interface as a function of time. The temperature gradient (G) applied is in the vertical direction. As opposed to the solidification in an undercooled pure metal, the liquid phase in this type of solidification is at a higher temperature than the dendrite in a solid phase.  Any numerically localized perturbation to the solid-liquid interface will grow as per the degree of constitutional undercooling. For the given temperature gradient $G$ and pulling velocity $\bm{\nu}_{p}$, the solid-liquid interface is unstable. This is evident from the interface $\phi=0$ evolution captured in the Subfigure~\ref{fig:Phi_SSC}. At t = 0, the solid seed starts to grow. By t = 35, the solid seed  has developed four primary branches due to the four-fold symmetry in the anisotropy considered. The growth of only three primary arms is captured in the simulation and the fourth arm pointing downward is not a part of the computation. These three primary arms have developed secondary arms of their own as seen in contour. By t = 40, the initial seed has grown into a larger size dendrite and several of its secondary arms are growing in the columnar fashion aligning with the direction in which constraints are applied, i.e, vertically. The phenomenon of secondary arms growing in the vertical direction becomes more obvious by t = 45 and t = 50 as seen in the Subfigure~\ref{fig:Phi_SSC}. The diffusion of the solute is governed by the same mechanism discussed in Section~\ref{sub:MultipleEquiaxedDendrite}. Like Subfigure~\ref{fig:Composition_ME}, in this case as well, the solute composition is maximum near the dendrite interface and it decreases as one moves away from the interface. This is clearly shown in the Subfigure~\ref{fig:Composition_SSC}.

\subsection{Growth of multiple columnar dendrites}
\begin{figure}[ht!]
  \centering 
  \subfloat[\large Composition\label{fig:Composition_MC}]{
       \includegraphics[width=0.45\textwidth]{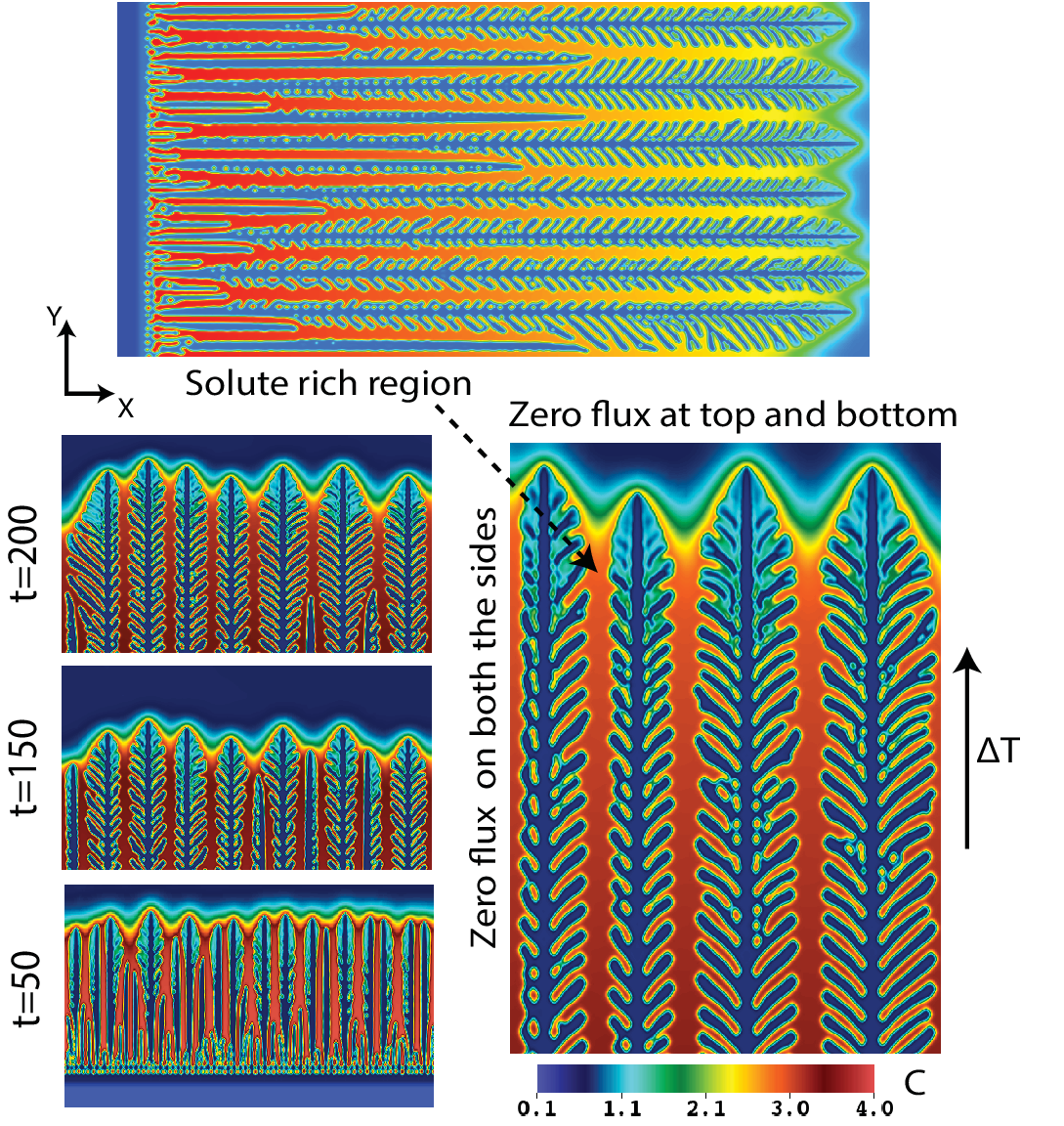}
  } 
   \subfloat[\large Zero level sets of the order parameter \label{fig:Phi_MC}]{
        \includegraphics[width=0.5\textwidth]{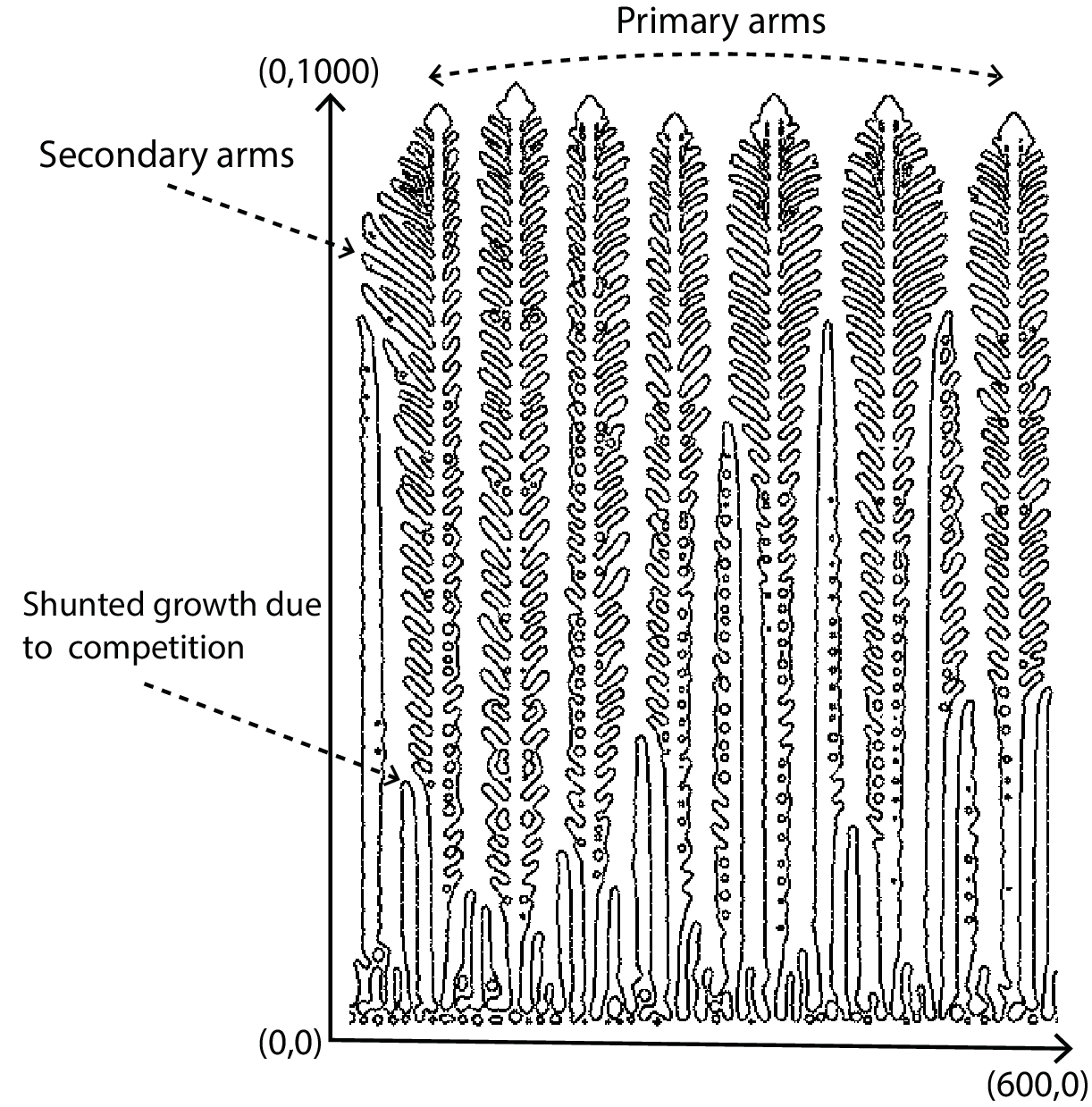}
  }
  \caption{Growth of multiple columnar dendrites. Shown are the (\ref{fig:Composition_MC}) Contours of the solute composition, $c$, at time t = 240, and (\ref{fig:Phi_MC}) Growth profile of multiple columnar dendrite interfaces at time t = 240. }
  \label{fig:ColumnarC_PHI}
\end{figure}

{\label{sub:MultiColumnar}}
Columnar dendrites growing out of a flat boundary are the next type of dendrites modeled. This type of dendrites occur in directional solidification processes which are of significance in manufacturing of alloy components by use of molds, where often the dendrite growth is primarily columnar.  The characteristic of the columnar dendrites are the growth of a colony of dendrites along the direction of the imposed temperature gradient and the pulling velocity. These dendrites then move at a steady-state tip velocity which is determined by the pulling velocity. The other geometric features of columnar dendrites which are typically observed in the experimental studies include primary arms, secondary arms, spacing between arms, growth competition between the neighboring dendrites, and the interaction of the diffusion field between neighboring dendrites. The growth of multi-columnar dendrites is governed by the same governing equations that were used to model the solidification of a binary alloy in the previous Section~\ref{sub:SingleEquiaxed}-\ref{sub:MultipleEquiaxedDendrite}, and so are the input parameters. However for the Initial Condition, we model an initial thin layer of solid solute on one of the boundaries, and set $\phi=1$ inside the initial solute layer thickness and $\phi=-1$ in the rest of the domain. Non-dimensional composition undercooling is set to $u_{c}=-1$ everywhere in the domain. For the Boundary Conditions, zero flux conditions are enforced on the boundaries for both the primal fields $u_{c}$ and $\phi$. Figure~\ref{fig:ColumnarC_PHI} depicts the results of multiple columnar dendritic growth. As can be seen, the phase-field contour interface ($\phi=0$) at time, t= 240, fills the entire domain with a colony of columnar dendrites. The initial planar interface with a solute composition seeds the dendrites that then compete to evolve into this columnar dendritic structure. The emergence of seven primary arms of the dendrites, as seen in Subfigure~\ref{fig:Phi_MC}, is the result of the competition for growth among several more initial primary arms of dendrites. Some of these underdeveloped arms can be seen in the lower portion of the figure. The primary arm spacing between these seven primary arms appears to converge to a constant value with the progress of the dendrites. Another feature of the columnar dendrites that we were able to capture is the growth of the secondary dendrite arms and their orientation. The concentration of the solute in the solid region is nearly uniform at a value of 0.56. The region surrounding these solid dendrites consists of the liquid melt where the highest concentration is about 4.0. Due to mass diffusion, there exists a gradient in the solute concentration in the liquid melt as can be seen in the Subfigure~\ref{fig:Composition_MC}. 

\subsection{Three dimensional growth of single equiaxed dendrite}
\begin{figure}[ht!]
  \centering
       \includegraphics[width=0.8\textwidth]{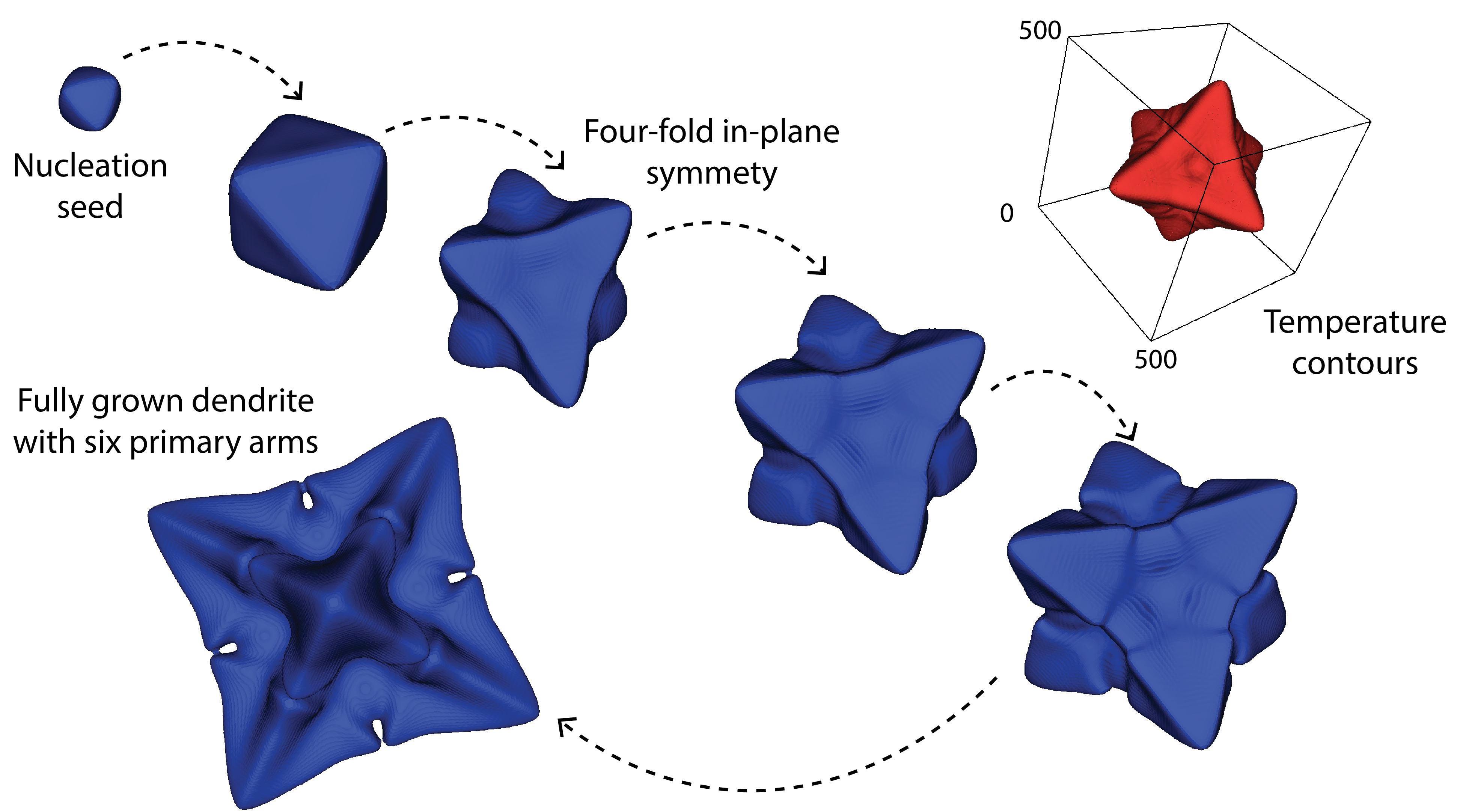}
  \caption{Evolution of a 3D single equiaxed dendrite. The dendritic structures shown in blue are the time evolution of the contours of the phase-field parameter, $\phi$, and the dendritic structure shown in red is the contour of the undercooling temperature, $u$, at one time instance.}
  \label{fig:3D_equiaxed}
\end{figure}

{\label{sub:3DSingleEquiaxed}}
We now simulate the phase-field model described by the Equations~\ref{eq:EquationU}-\ref{eq:EquationPHI} in three dimensions. The initial solidification conditions consist of a small spherical seed placed in a $500\times500\times500$ computational domain with an adaptive mesh and the smallest element length measure of $\Delta x=2.0$ and uniform time step  $\Delta t=0.12$.  The seed is placed in an undercooled melt of pure metal. The input phase-field parameters used in the simulation are listed in Table~\ref{table:PureMeltParameters}. The evolution of a spherical seed into a four-fold symmetric dendrite structure is shown in the Figure~\ref{fig:3D_equiaxed}. The fully grown dendrite in blue at the end of the evolution has six primary arms as dictated by the surface anisotropy. The distribution of an undercooling temperature $u$ on a fully grown dendrite, shown in red, can also be seen in the figure. Three dimensional implementation of phase-field models of dendrites is computationally very demanding, due to the need to resolve the fine scale primary-arm and secondary-arm structures. To make this problem computationally tractable, we use local mesh adaptivity (h-adaptivity), implicit time-stepping schemes with adaptive time-step control, and domain decomposition using MPI.  

\section{Conclusion}{\label{sub:DiscussionConclusion}}

\noindent In this work, we present a detailed development of the requisite numerical aspects of dendritic solidification theory, and the computational implementation of the corresponding phase-field formulations to model dendritic growth in pure metals and binary alloys. A wide variety of physically relevant dendritic solidification patterns are modeled by solving the governing equations under various initial conditions and boundary conditions. To validate the numerical framework, we simulate the classical four-fold symmetric dendrite shape occurring in undercooled pure melt and compare the numerically computed dendritic tip velocities with the corresponding analytical values obtained using a Green's function approach. Further, this problem is used as a basis for performing error convergence studies of dendritic tip velocity and dendritic morphology (primal fields of temperature and order parameter). Further, using this numerical framework, various types of dendritic solidification patterns like multi-equiaxed, single columnar and multi-columnar dendrites are modeled in two-dimensional and three-dimensional computational domains.  

The distinguishing aspects of this work are - a unified treatment of both pure-metals and alloys; novel numerical error estimates of dendritic tip velocity; and the study of error convergence of the primal fields of temperature and the order parameter with respect to the numerical discretization. To the best of our knowledge, this is a first of its kind study of numerical convergence of the phase-field equations of dendritic growth in a finite element method setting, and a unified computational framework for modeling a variety of dendritic structures relevant to solidification in metallic alloys. 
  
\section*{Acknowledgement}
The authors would like to thank Prof. Dan Thoma (University of Wisconsin-Madison) and Dr. Kaila Bertsch (University of Wisconsin-Madison; now at Lawrence Livermore National Laboratory) for very useful discussions on dendritic growth and microstructure evolution in the context of additive manufacturing of metallic alloys. 

\section*{Declarations}    
\textbf{Conflict of interest}: The authors declare that they have no conflict of interest.


\bibliography{main}

\end{document}